\documentclass[12pt,a4paper]{article}

\usepackage{amsmath}
\usepackage{amsthm}
\usepackage{amsfonts}
\usepackage{amssymb}
\usepackage{graphicx}
\usepackage[margin=2.6cm]{geometry}
\usepackage{enumitem}
\usepackage[T1]{fontenc}
\usepackage[utf8]{inputenc}
\usepackage{tikz}
\usetikzlibrary{calc}
\usepackage{courier}
\usepackage{xcolor}

\def\R{{\mathbb R}}
\def\E{{\mathbb E}}

\def\F{{\mathcal F}}

\def\1{{\bf 1}}
\def\0{{\bf 0}}

\def\aff{\operatorname{aff}}
\def\conv{\operatorname{conv}}
\def\cone{\operatorname{cone}}
\def\lin{\operatorname{lin}}
\def\inte{\operatorname{int}}
\def\relint{\operatorname{relint}}
\def\rk{\operatorname{rk}}

\def\se{\subseteq}
\def\imp{\Rightarrow}

\providecommand{\texorpdfstring}[2]{#1}

\def\pole#1{\fill[white] (#1) circle (4.2pt);\fill[blue!70!black] (#1) circle (3pt);}
\def\gen#1{\fill[white] (#1) circle (3.8pt);\fill[black!65] (#1) circle (2.6pt);}
\def\genw#1{\fill[white] (#1) circle (3.8pt);\fill[black!35] (#1) circle (2.6pt);}
\def\rightangle#1#2#3{%
  \coordinate (raA) at ($(#1)!0.22!(#2)$);%
  \coordinate (raB) at ($(#1)!0.22!(#3)$);%
  \draw[thin] (raA) -- ($(raA)+(raB)-(#1)$) -- (raB);}

\newtheorem{theo}{Theorem}[section]
\newtheorem{coro}[theo]{Corollary}
\newtheorem{lm}[theo]{Lemma}
\newtheorem{prop}[theo]{Proposition}

\newtheorem{rem}[theo]{Remark}
\newtheorem{de}[theo]{Definition}

\newenvironment{proo}{\begin{proof}}{\end{proof}}

\begin{document}

\begin{center}
\LARGE{\bf Thales's theorem in higher dimensions}\\
\vspace{0.5 cm}
\normalsize{Jan Pavl\'ik} \\
\vspace{0.5 cm}
\small Faculty of Mechanical Engineering, Brno University of Technology,\\
Technick\'a 2896/2, 616 69 Brno, Czech Republic\\
{\tt pavlik@fme.vutbr.cz}\\
\end{center}

\begin{abstract}
In the Euclidean plane, Thales's theorem is equivalent to the statement that a
quadrilateral with two opposite right angles has a circumscribed circle. We generalize
the latter to higher dimensions: a convex polytope with two distinguished vertices, the
\emph{poles}, such that the faces spanned by the two poles are perpendicular, is
inscribed in the sphere whose diameter joins its poles. These \emph{orthogonally bipolar}
polytopes are exactly the bipolar polytopes of a companion paper that arise as cutouts
by orthogonal projections; in particular, the interval polytope of $Q$ has such an
inscribed realization exactly when $Q$ has a perfectly centered realization. As an application we show that for a polytope $P$ and a point
$B$, the balls with diameters $[V,B]$, $V$ a vertex of $P$, cover $P$, and that the power
diagram of their boundary spheres refines this covering to a tiling. Further, we
investigate the properties of orthogonally bipolar polytopes, their characterization,
and their application to tilings: the cutouts of $P$ by the orthogonal projections of
$B$ onto its faces tile $P$ by orthogonally bipolar polytopes exactly when every such
projection lies in the relative interior of its face, i.e.\ when $P$ is perfectly
centered about $B$. We also observe that the Thales condition underlies Bessel's inequality
and Parseval's identity for orthogonal expansions, and characterize when the best approximations of a vector by the
subfamilies of a basis form a Hjelmslev cuboid.
\end{abstract}


\section{Preliminaries}\label{sec:prelim}

\subsection{Motivation}

Thales's theorem is among the oldest results of geometry: in the Euclidean plane, the
angle subtended by a diameter of a circle at a point of the circle is a right angle. It
ties together two basic notions, orthogonality and the circle, and its converse
characterizes the points of a circle by right angles. In higher dimensions orthogonality
is a richer notion: not only lines, but whole faces of a polytope can be perpendicular.
This raises a natural question: which configurations of perpendicular faces force a
polytope to be inscribed in a sphere? The plane suggests the setting through another
statement: a quadrilateral whose two opposite angles are
right angles has a circumscribed circle. The two statements are equivalent,
in the sense that either can be proved directly from the other. We show
that the latter has a multidimensional counterpart. One of its consequences
is that the two vertices at which the quadrilateral does \emph{not} have a
right angle realize the diameter of the quadrilateral. We then
show that a generalized version of the multidimensional statement has an
analogous consequence.

The Thales condition has a familiar counterpart in analysis. If $P$ is the orthogonal
projection onto a closed subspace of a Hilbert space and $\nu$ is a vector, then
$\langle P\nu,\nu-P\nu\rangle=0$: all orthogonal projections of $\nu$ lie on the sphere with
diameter $[0,\nu]$. For instance, the partial sums $S_kf$ of the Fourier series of
$f\in L^2$ lie on the sphere with diameter $[0,f]$ and converge to its antipode $f$; the
identity $\|S_kf\|^2=\langle S_kf,f\rangle$ behind Bessel's inequality is exactly the Thales
condition. Projecting $\nu$ onto the spans of all subfamilies of an orthonormal basis gives
the vertices of a box with opposite vertices $0$ and $\nu$. For a non-orthogonal basis the
projections still lie on the Thales sphere, and Proposition~\ref{prop:bases} shows that
they form a Hjelmslev cuboid, in general not a box, exactly when the coefficients of the
best approximations never change their signs. The polytopes studied here can thus be seen
as the finite, combinatorially structured part of a geometry of orthogonal projections.

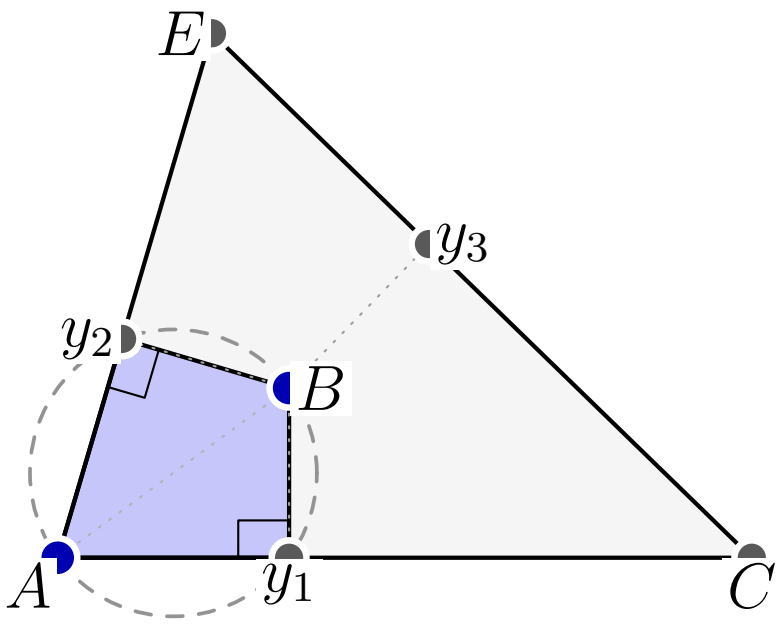
\begin{figure}[htbp]
\centering
\begin{tikzpicture}[scale=1.02,line cap=round,line join=round]
  \coordinate (A)  at (0,0);
  \coordinate (C)  at (4.5,0);
  \coordinate (E)  at (1,3.4);
  \coordinate (B)  at (1.5,1.1);
  \coordinate (y1) at (1.5,0);
  \coordinate (y2) at (0.4172,1.4185);
  \coordinate (y3) at (2.4068,2.0334);
  \fill[gray!8] (A)--(C)--(E)--cycle;
  \fill[blue!30,fill opacity=0.7] (A)--(y1)--(B)--(y2)--cycle;
  \draw[thick] (A)--(C)--(E)--cycle;
  \draw[thick] (A)--(y1)--(B)--(y2)--cycle;
  \draw[gray!75,dotted] (B)--(y1) (B)--(y2) (B)--(y3);
  \draw[gray!85,dashed,line width=0.7pt] (0.75,0.55) circle (0.9301);
  \draw[gray!60,dotted] (A)--(B);
  \rightangle{y1}{A}{B}
  \rightangle{y2}{A}{B}
  \pole{A} \pole{B} \gen{y1} \gen{y2} \gen{y3} \gen{C} \gen{E}
  \node[below left, fill=white,inner sep=1pt] at (A)  {$A$};
  \node[right,      fill=white,inner sep=1pt] at (B)  {$B$};
  \node[below,      fill=white,inner sep=1pt] at (y1) {$y_1$};
  \node[left,       fill=white,inner sep=1pt] at (y2) {$y_2$};
  \node[right,      fill=white,inner sep=1pt] at (y3) {$y_3$};
  \node[below,      fill=white,inner sep=1pt] at (C)  {$C$};
  \node[left,       fill=white,inner sep=1pt] at (E)  {$E$};
\end{tikzpicture}
\caption{The planar case, and the reference picture for everything that
follows. For the triangle $ACE$ and the point $B$ the points
$y_1,y_2,y_3$ are the orthogonal projections of $B$ onto the three edges;
$y_1,y_2$ lie on the faces through $A$. The quadrilateral $A,y_1,B,y_2$ is
the cutout of the orthogonal configuration of $B$ at $A$: it has right
angles at the generators $y_1,y_2$, its poles are $A$ and $B$, and
the circle with diameter $AB$ is circumscribed to it, as Proposition
\ref{hjel} predicts. Here all three feet fall into the relative interiors
of the edges, so the projections form a configuration; a case where they
do not is drawn in Figure~\ref{fig:cover}.}
\label{fig:hjelcut}
\end{figure}

A second motivation, besides this geometric question, comes from combinatorics. In a
companion paper \cite{Bipolar} we study \emph{bipolar polytopes}, that is, polytopes with
two distinguished vertices, the \emph{poles} $\sigma,\nu$, such that every
facet contains exactly one of them and a certain intersection property
holds; the vertex set then carries a natural lattice order whose intervals
are precisely the faces. Such polytopes arise geometrically as
\emph{cutouts} of the apex of a pyramid: taking the apex and an interior
point $B$ as the two poles and the orthogonal projections $x_G$ of $B$ onto
the faces $G$ at the apex as the remaining vertices, one obtains, under a
condition on the normal fan and the face fan of the base
\cite[Definition~5.1]{Bipolar}, a convex polytope whose vertex lattice is
the face lattice of the base Proposition~\ref{prop:cutout}. Orthogonality is
the ingredient of this construction that we study here: such a cutout is an
orthogonally bipolar polytope, and the sphere with diameter $[v,B]$ is
circumscribed to it by Theorem~\ref{theo:obpthales}. Such a realization
exists exactly when the combinatorial type has a perfectly centered realization, by a
theorem of Broadie \cite{Broadie} (see Theorem~\ref{conj:I}); this is the case in
dimension at most three \cite[Theorem~2]{GrunbaumAnti}, \cite[Theorem~3.4]{Dobbins}, but not for every
$4$-polytope \cite[Theorem~3.9]{Dobbins}. The present paper studies
the orthogonal structure behind that question; the
Hjelmslev quadrilateral of Section~\ref{sec:thales} is its planar case and
the Hjelmslev $n$-cuboids its Boolean case
(Figure~\ref{fig:hjelcut}). Taking the cutouts at all vertices of a polytope at once
leads to tilings by orthogonally bipolar polytopes; Section~\ref{sec:pedal} shows, using Broadie's theorem,
that they exist exactly for the points of the orthogonal kernel
(Theorem~\ref{thm:pedal}), i.e.\ the points about which the polytope is perfectly
centered, and that they are then the cells of the power diagram of the
Thales spheres.

\medskip
\noindent\textbf{Main results.}
\begin{description}[leftmargin=1em,style=nextline]
\item[Thales's theorem (Theorem~\ref{theo:obpthales}).]
Every orthogonally bipolar polytope is inscribed in the sphere whose diameter joins its
poles.

\item[Orthogonal cutouts (Propositions~\ref{prop:cut-obp} and~\ref{prop:obp-cut}).]
The orthogonally bipolar polytopes are exactly the admissible cutouts formed by the
orthogonal projections of a point onto the faces through a vertex.

\item[Inscribed interval polytopes (Theorem~\ref{prop:obpreal}).]
The interval polytope of a polytope $Q$, the bipolar polytope $W(\F(Q))$ of
\cite{Bipolar}, has an orthogonally bipolar realization, inscribed in a sphere, if and
only if the combinatorial type of $Q$ has a perfectly centered realization; every such
realization is the cutout of the apex of a pyramid over $Q$ by orthogonal projections.
Since perfect centering is inherited by the interval polytope (Lawrence), all iterated
interval polytopes then have such realizations (Corollary~\ref{cor:iterate}).

\item[Covering (Theorem~\ref{thm:cover} and Corollary~\ref{cor:power}).]
For a polytope $P$ and a point $B$, the balls with diameters $[V,B]$, $V$ a vertex of
$P$, cover $P$; the power diagram of their boundary spheres tiles $P$ by cells lying in
these balls.

\item[Orthogonal tilings (Theorem~\ref{thm:pedal}).]
The orthogonal cutouts of a polytope $P$ at a point $B$ are orthogonally bipolar and
tile $P$ face-to-face if and only if $B$ lies in the orthogonal kernel of $P$, i.e.\ $P$ is
perfectly centered about $B$.
\end{description}

\subsection{General setting}\label{sec:notation}

This subsection fixes notation and recalls the notions used throughout the paper.
Standard notions are recalled briefly, with references; the notions and results taken
from the companion paper \cite{Bipolar} are stated precisely, in particular the facts
(F1)--(F3) used in the proofs; and the notions specific to this paper that are used in
several sections (cutout configurations, orthogonal configurations, Thales spheres) are
introduced here as well.

\paragraph{Affine space.} We work in the $n$-dimensional real affine space. Fixing an
origin $O$, we identify it with $\R^n$ and a point $x$ with its position vector $x-O$.
The difference $x-y$ of two points is a vector, and an affine combination
$\sum\lambda_ix_i$ with $\sum\lambda_i=1$ is a point that does not depend on the
choice of $O$. The vectors carry the standard inner product
$\langle\cdot,\cdot\rangle$ and the norm $\|x\|=\langle x,x\rangle^{1/2}$, and the
points the metric $\delta(x,y)=\|x-y\|$; thus our space is the Euclidean space
$\R^n$.\footnote{Our notation and conventions agree with those of the companion paper \cite{Bipolar}.}
For the standard notions of polyhedral combinatorics we follow
\cite{Carr2005,Schrijver95,Ziegler,Grunbaum}. An \emph{affine subspace} is a set
$\{x\in\R^n:Ax=b\}$ with a matrix $A$ and a vector $b$; equivalently, a set closed under
affine combinations. A (closed)
\emph{halfspace} is a set $\{x\in\R^n:\langle\alpha,x\rangle\le b\}$ with
$\alpha\in\R^n\setminus\{O\}$ and $b\in\R$. A set is \emph{convex} if with any two
points $x,y$ it contains the segment $\{tx+(1-t)y:t\in[0,1]\}$, i.e.\ if it is closed
under convex combinations.
For $S\subseteq\R^n$ let $\aff S$ be the \emph{affine hull}, the smallest
affine subspace containing $S$, i.e.\ the set of all affine combinations
$\sum\lambda_ix_i$ with $x_i\in S$, $\lambda_i\in\R$, $\sum\lambda_i=1$
(negative coefficients are allowed, in contrast to $\conv S$);
$\conv S$ is the convex hull, $\cone S=\{\sum\lambda_ix_i:\ x_i\in S,\
\lambda_i\ge0\}$ the cone over $S$, and $\lin S$ the linear span. We write
$\relint S$ and $\inte S$ for the relative interior (the interior within
$\aff S$) and the interior, and put $\dim S=\dim\aff S$. The closed ball
and the sphere (the boundary of the ball) of radius $r$ with centre $X$ are
denoted by $B(X,r)$ and $S(X,r)$.

\paragraph{Polytopes.} A \emph{convex polyhedron} is an intersection of finitely many
halfspaces, and a \emph{polytope} is a bounded convex polyhedron, equivalently the
convex hull of finitely many points. The dimension of a polytope $P$ is $\dim\aff P$,
i.e.\ the least $d$ such that $P$ lies in an affine subspace $\{x:Ax=b\}$ with $A$ of
rank $n-d$. A \emph{face} of $P$ is $P$ itself or the intersection of $P$ with a
\emph{supporting hyperplane} $\{x:\langle\alpha,x\rangle=b\}$, i.e.\ one such that
$P\subseteq\{x:\langle\alpha,x\rangle\le b\}$; the empty set is a face, and every face
is again a polytope. A \emph{$k$-face} is a face of dimension $k$; the $0$-faces,
$1$-faces and $(\dim P-1)$-faces are the \emph{vertices}, \emph{edges} and
\emph{facets} of $P$, and $V(P)$ denotes the vertex set.

\paragraph{Perpendicularity.} A \emph{vertex vector} of a polytope $W$ is
a difference $x-y$ of two of its vertices. Two faces $F,G$ of $W$ are
\emph{perpendicular}, written $F\perp G$, if every vertex vector of $F$ is
orthogonal to every vertex vector of $G$; equivalently, the direction
spaces of $\aff F$ and $\aff G$ are orthogonal. This is the relation used
in the definition of a Hjelmslev cuboid.

\paragraph{Direction space.} For a face $F$ of a polytope we denote by
$U_F$ the direction space of $\aff F$, i.e.\ $U_F=\{x-y:x,y\in\aff F\}$;
equivalently, $U_F$ is the linear span of the directions of the edges of
$F$ at any vertex of $F$.

\paragraph{Poles.} Let $W$ be a convex polytope with two distinct vertices
$\sigma\ne\nu$, the \emph{poles} (\emph{south} $\sigma$, \emph{north}
$\nu$). A face is a \emph{south} face if it contains $\sigma$ and a
\emph{north} face if it contains $\nu$; for a vertex $x$ we write
$[\sigma,x]$ and $[x,\nu]$ for the faces spanned by $\sigma,x$, resp.\
$x,\nu$.

\paragraph{Face lattices.} For a polytope $W$ let $\F(W)$ be its
\emph{face lattice}: the set of all faces of $W$, including $\emptyset$ and
$W$, ordered by inclusion; it is a lattice. The \emph{rank} of a face is
$\rk F=\dim F$, with $\rk\emptyset=-1$; thus vertices have rank $0$ and
facets rank $n-1$. With this rank function $\F(W)$ is an abstract
polytope of rank $n$ \cite[\S2A]{McMullenSchulte}. Two polytopes are \emph{combinatorially equivalent}, or
of the same \emph{combinatorial type}, if their face lattices are
isomorphic, and a polytope $Q'$ with the same face lattice as $Q$ is a
\emph{realization} of the type of $Q$. In a finite lattice with least
element $0$ and greatest element $1$ the \emph{atoms} are the elements
covering $0$ and the \emph{coatoms} the elements covered by $1$.

\paragraph{Vertex figures.} For a vertex $v$ of an $n$-polytope $W$ the faces
containing $v$ form a lattice isomorphic to the face lattice $\F(W/v)$ of the
\emph{vertex figure} $W/v$ of $W$ at $v$; under this isomorphism $v$ corresponds to
$\emptyset$ and $W$ to $W/v$.

\paragraph{Bipolar polytopes: facts from \cite{Bipolar}.} Let $W$ be a convex polytope
with poles $\sigma\ne\nu$. Then $W$ is a \emph{geometric bipolar polytope} with poles
$\sigma,\nu$ if every facet of $W$ contains exactly one of the poles, and for every
south face $H$ the family of north faces meeting $H$ is closed under intersection, and
dually \cite[Definitions~3.1 and~3.2]{Bipolar}; both conditions reappear below as (O1),
(O2) of Definition~\ref{de:obp}. We use the following three results of
\cite{Bipolar}, for a geometric bipolar polytope $W$ with poles $\sigma,\nu$.
\begin{enumerate}[label=(F\arabic*)]
\item \cite[Theorem~3.11]{Bipolar} Write $x\le y$ for vertices $x,y$ if $y$ lies in
      every north face containing $x$. Then $(V(W),\le)$ is a lattice, the
      \emph{vertex lattice} of $W$, with least element $\sigma$ and greatest element
      $\nu$. Its intervals $[x,y]$, $x\le y$, are exactly the vertex sets of the
      nonempty faces of $W$, and the facets of $W$ are the faces $[a,\nu]$, $a$ an
      atom, and $[\sigma,c]$, $c$ a coatom.
\item \cite[Theorem~3.12]{Bipolar} For $x<y$ the face $[x,y]$ is a geometric bipolar
      polytope with poles $x,y$, and its vertex lattice is the interval $[x,y]$ of
      $(V(W),\le)$.
\item \cite[Proposition~3.13]{Bipolar} The map $y\mapsto[\sigma,y]$ is an isomorphism
      of $(V(W),\le)$ onto the lattice of faces containing $\sigma$, i.e.\
      $(V(W),\le)\cong\F(W/\sigma)$. In particular the atoms of $(V(W),\le)$ are the
      vertices joined to $\sigma$ by an edge.
\end{enumerate}

\paragraph{Cutout configurations and cutouts.} Let $P\se\R^n$ be an
$n$-polytope, $n\ge2$, and $v$ a vertex of $P$. A \emph{cutout
configuration at $v$} consists of a point $B\in\inte P$ and points
$x_G\in\relint G$ for the faces $G\ni v$ with $\{v\}\ne G\ne P$; put
$x_{\{v\}}=v$ and $x_P=B$. The \emph{$(B,v)$-cutout} of $P$ is
\[
C=\conv\{x_G:\ v\in G\in\F(P)\}.
\]
It is \emph{admissible} if the vertices of $C$ are exactly the (pairwise
distinct) points $x_G$ and $C$ is a geometric bipolar polytope with poles
$v$ and $B$. If $P$ is a pyramid $\conv(Q'\cup\{v\})$ over a polytope $Q'$ with
apex $v\notin\aff Q'$, then $C$ is the \emph{apex cutout} of the pyramid. (A cutout configuration at $v$ is a
configuration at $v$ in the sense of \cite[Definition~5.7]{Bipolar}.)

\paragraph{Orthogonal configurations.} Given
$B\in\inte P$ and a vertex $v$ of $P$, let $x_G$ be the \emph{orthogonal
projection} of $B$ onto $\aff G$, for every face $G\ni v$. Then
$x_{\{v\}}=v$ and $x_P=B$ automatically. If $x_G\in\relint G$ for all these
faces, $(x_G)$ is a cutout configuration at $v$, called the
\emph{orthogonal configuration} of $B$ at $v$; in general, however, the
projections need not lie in the relative interiors of the faces
(Figure~\ref{fig:cover}). In all cases $x_G-B\perp U_E$
for faces $E\se G$ (\emph{orthogonal coplanarity}: $x_G-B$ is orthogonal
to $U_G\supseteq U_E$); in particular all $x_G$ with $G\supseteq e$ lie
in the hyperplane $B+U_e^{\perp}$ orthogonal to an edge $e$ at $v$.

\paragraph{Tilings.} A finite family of $n$-polytopes $C_1,\dots,C_k$ \emph{tiles} an
$n$-polytope $P$, or forms a \emph{tiling} of $P$, if $C_1\cup\dots\cup C_k=P$ and the
interiors of the $C_i$ are pairwise disjoint.

\paragraph{Voronoi and power diagrams.} For a finite set $S\subset\R^n$ of sites, the
\emph{Voronoi cell} of $s\in S$ consists of the points at least as close to $s$ as to
any other site, and these cells form the \emph{Voronoi diagram} of $S$; its dual, formed
by the convex hulls of the sites whose cells share a vertex, is the \emph{Delaunay
subdivision}, a triangulation for sites in general position. If each site carries a
weight $w_s$, the \emph{power} of a point $p$ with respect to $s$ is
\[
\operatorname{pow}_s(p)=\|p-s\|^2-w_s ;
\]
for $w_s=r^2>0$ it is the power of $p$ with respect to the sphere with centre $s$ and
radius $r$. Replacing the distance by the power gives the \emph{power diagram}, whose
cells are separated by the radical hyperplanes of the spheres; with equal weights it
is the Voronoi diagram \cite{Aurenhammer87,Aurenhammer91}. Power diagrams of Thales spheres appear in
Corollary~\ref{cor:power}.

\paragraph{Thales spheres and balls.} For distinct points $X,Y\in\R^n$ the
\emph{sphere with diameter $[X,Y]$} is
$S\bigl(\tfrac12(X+Y),\tfrac12\delta(X,Y)\bigr)$, the $(n-1)$-sphere with
centre the midpoint of $[X,Y]$ on which $X$ and $Y$ are antipodal; the
\emph{closed ball with diameter $[X,Y]$} is
\[
D_{XY}=\{Z\in\R^n:\ \langle Z-X,\,Z-Y\rangle\le0\}.
\]
Indeed, with the midpoint $M=\tfrac12(X+Y)$ and $\rho=\tfrac12\delta(X,Y)$ we have
$Z-X=(Z-M)+(M-X)$, $Z-Y=(Z-M)-(M-X)$ and $\|M-X\|=\rho$, hence
\[
\langle Z-X,\,Z-Y\rangle=\|Z-M\|^2-\rho^2 ;
\]
so $D_{XY}$ is the closed ball $B(M,\rho)$, and the sphere with diameter $[X,Y]$ is the
set where $\langle Z-X,Z-Y\rangle=0$. In particular $Z\ne X,Y$ lies on the sphere iff
$\angle XZY=\tfrac\pi2$ (the $n$-dimensional Thales condition). The same formula
defines $D_{XY}$ also for $X=Y$, where it is the single point $X$. A sphere is \emph{circumscribed} to a
polytope if it contains all its vertices.

\section{Thales's theorem}\label{sec:thales}

The theorem is well known and there are several ways to prove it; we only
recall it together with its converse.

\begin{prop}\label{thal1}
(Thales's theorem -- the direct version) Let $A,B,C\in\E_2$ be points on a
circle of radius $r$ such that $\delta(A,C)=2r$, and $B\ne A,C$. Then the angle
$\angle(ABC)$ is $\frac{\pi}{2}$.
\end{prop}

\begin{prop}\label{thal2}
(Thales's theorem -- the converse version) The circumcentre of a right
triangle lies at the midpoint of its hypotenuse.
\end{prop}

These propositions are well known, and each can be proved as a consequence
of the other.

E.g., the direct version follows from the converse as follows. Let
$A,B,C\in\E_2$ lie on a circle with centre $M$, with $AC$ a diameter and
$B\ne A,C$. Let $B'$ be the foot of the perpendicular from $C$ to the line
$AB$. If $B'=A$, then $CA\perp AB$, so the line $AB$ is tangent to the
circle at $A$ and cannot meet it again at $B$; hence $B'\ne A$, and the
triangle $\triangle AB'C$ has a right angle at $B'$. By
Proposition~\ref{thal2} its circumcentre is the midpoint $M$ of $AC$, so
$B'$ lies on the circle. The line $AB$ meets the circle exactly in $A$ and
$B$, and $B'\ne A$ lies on both; hence $B'=B$ and
$\angle(ABC)=\angle(AB'C)=\frac\pi2$.

Conversely, suppose the direct version holds and let $\triangle ABC$ have
a right angle at the vertex $B$, with $M$ the midpoint of the hypotenuse
$AC$. The ray $\overrightarrow{MB}$ meets the circumcircle of
$\triangle ABC$ at a point $D$. By the direct version the angle
$\angle(ADC)$ is a right angle; since $B$ and $D$ both lie on the ray
$\overrightarrow{MB}$ and a point further along the ray sees the segment
$AC$ under a smaller angle, we get $D=B$. Hence $B$ lies on the circle and
the circle is circumscribed to the triangle.

\begin{de}
A \emph{Hjelmslev quadrilateral} is a polygon with vertices
$A,B,C,D\in\E_2$ in the cyclic order given alphabetically such that the
angles $\angle(ABC)$ and $\angle(ADC)$ are right angles. The points $A$
and $C$ are called the \emph{poles} and $B,D$ the \emph{generators}.
\end{de}

\begin{prop}\label{hjel}
Let $Q$ be a Hjelmslev quadrilateral with vertices $A,B,C,D$ and poles
$A,C$. Then $Q$ has a circumscribed circle whose centre is the midpoint of
$AC$.
\end{prop}

We show that this statement and the converse Thales's theorem imply each
other.

\begin{prop}
Propositions~\ref{thal2} and~\ref{hjel} can be proved from each other.
\end{prop}
\begin{proo}
``Proposition~\ref{thal2} $\imp$ Proposition~\ref{hjel}'': assume
Proposition~\ref{thal2} and let $Q$ be a Hjelmslev quadrilateral with
vertices $A,B,C,D$ and poles $A,C$; let $M$ be their midpoint. By
Proposition~\ref{thal2}, $r=\delta(A,M)$ is a radius of the circumcircles
of both $\triangle ABC$ and $\triangle ADC$; hence the circle $S(M,r)$ is
a circumcircle of the quadrilateral $A,B,C,D$.

``Proposition~\ref{hjel} $\imp$ Proposition~\ref{thal2}'': let
$\triangle ABC$ be a right triangle with $M$ the midpoint of its
hypotenuse. The point $D=M+(M-B)$ completes a Hjelmslev quadrilateral
$A,B,C,D$, which by Proposition~\ref{hjel} has a circumcircle with centre
$M$; that circle is also a circumcircle of $\triangle ABC$.
\end{proo}

Thus the statement of Proposition~\ref{hjel} is equivalent to Thales's
theorem. We now look for its multidimensional version.

\section{Orthogonally bipolar polytopes}\label{sec:obp}

\begin{de}[orthogonally bipolar polytope]\label{de:obp}
Let $W\subseteq\R^n$ be a convex polytope with two distinct vertices
$\sigma\ne\nu$, called the \emph{poles} (\emph{south} $\sigma$,
\emph{north} $\nu$). Then $W$ is an \emph{orthogonally bipolar polytope} with poles $\sigma,\nu$ if
\begin{enumerate}[label=(O\arabic*)]
\item every facet of $W$ contains exactly one of the poles;
\item for every south face $H$ the family of north faces meeting $H$ is closed
      under intersection, and dually (\S\ref{sec:notation});
\item $[\sigma,a]\perp[a,\nu]$ for every atom $a$ of the vertex lattice
      $(V(W),\le)$.
\end{enumerate}
An orthogonally bipolar polytope whose vertex lattice is the Boolean lattice $B_n$ of the unit
$n$-cube is called a \emph{Hjelmslev $n$-cuboid} ($n$-HC); its poles are
the vertices $\0$ and $\1$. See \S\ref{sec:boolean} for a characterization
of that case.
\end{de}

\begin{rem}\label{rem:obp-prep}
Conditions (O1) and (O2) say that $W$ is a geometric bipolar polytope with
poles $\sigma,\nu$ (\S\ref{sec:notation}). By (F1) the vertices then carry the
vertex lattice $(V(W),\le)$, whose intervals are exactly the faces of $W$; the
facets are the intervals $[a,\nu]$ with $a$ an atom and $[\sigma,c]$ with $c$ a
coatom, and $(V(W),\le)\cong\F(W/\sigma)$ by (F3). In particular the atoms of $(V(W),\le)$
correspond to the edges of $W$ at $\sigma$, so (O3) is a condition on
these edges.

Condition (O3) is stated for atoms only; the coatom version
$[\sigma,c]\perp[c,\nu]$ for a coatom $c$ is a consequence, because
$\perp$ passes to subfaces: if $a\le c$ is an atom, then
$[c,\nu]\subseteq[a,\nu]$, so $U_{[c,\nu]}\subseteq U_{[a,\nu]}$, and
$[\sigma,a]\perp[a,\nu]$ gives $[\sigma,c]\perp[c,\nu]$ after taking the
span of the edge directions at $\sigma$.
\end{rem}

\begin{figure}[htbp]
\centering
\begin{tikzpicture}[line join=round,line cap=round]
  \begin{scope}[scale=1.10]
    \coordinate (s)  at (0,0);
    \coordinate (e1) at (2,0);
    \coordinate (n)  at (2,2);
    \coordinate (e2) at (0,2);
    \draw[gray,dashed] (1,1) circle (1.41421);
    \draw[gray,densely dashed] (s)--(n);
    \draw[thick] (s)--(e1)--(n)--(e2)--cycle;
    \rightangle{e1}{s}{n}
    \rightangle{e2}{s}{n}
    \pole{s} \pole{n} \gen{e1} \gen{e2}
    \node[below left=3pt,inner sep=0pt]  at (s)  {$\sigma$};
    \node[above right=3pt,inner sep=0pt] at (n)  {$\nu$};
    \node[below right=3pt,inner sep=0pt] at (e1) {$\epsilon_1$};
    \node[above left=3pt,inner sep=0pt]  at (e2) {$\epsilon_2$};
  \end{scope}
  \begin{scope}[xshift=5.2cm,yshift=-0.35cm,scale=1.45]
    \coordinate (S)   at (0,0);
    \coordinate (A)   at (1.6,0);
    \coordinate (c)   at (1.6,1.6);
    \coordinate (f01) at (0,1.6);
    \coordinate (hid) at (0.8,0.56);
    \coordinate (b10) at (2.4,0.56);
    \coordinate (n)   at (2.4,2.16);
    \coordinate (b01) at (0.8,2.16);
    \fill[orange!12] (S)--(A)--(c)--(f01)--cycle;
    \fill[blue!40] (A)--(c)--(n)--(b10)--cycle;
    \draw[gray,dashed] (S)--(hid) (hid)--(b10) (hid)--(b01);
    \draw[thick] (S)--(A)--(c)--(f01)--cycle;
    \draw[thick] (A)--(b10)--(n)--(c);
    \draw[thick] (n)--(b01)--(f01);
    \draw[gray,densely dashed] (S)--(n);
    \draw[blue!70!black,very thick] (S)--(A);
    \draw[orange!80!black,very thick] (c)--(n);
    \draw[blue!70!black,line width=0.5pt] (1.152,0.000) -- (1.152,0.448) -- (1.600,0.448);
    \draw[blue!70!black,line width=0.5pt] (1.152,0.000) -- (1.376,0.157) -- (1.824,0.157);
    \draw[orange!80!black,line width=0.5pt] (1.824,1.757) -- (1.824,1.309) -- (1.600,1.152);
    \draw[orange!80!black,line width=0.5pt] (1.824,1.757) -- (1.376,1.757) -- (1.152,1.600);
    \pole{S} \pole{n}
    \gen{A} \gen{c} \gen{f01} \gen{b10} \gen{b01} \genw{hid}
    \node[below left=3pt,inner sep=0pt]  at (S)   {$\sigma$};
    \node[above right=3pt,inner sep=0pt] at (n)   {$\nu$};
    \node[below=4pt,inner sep=0pt]       at (A)   {$a=\epsilon_1$};
    \node[below left=5pt,inner sep=0pt]  at (c)   {$c$};
    \node[left=4pt,inner sep=0pt]        at (f01) {$\epsilon_2$};
    \node[above left=2pt,inner sep=0pt,gray] at (hid) {$\epsilon_3$};
  \end{scope}
\end{tikzpicture}
\caption{Two orthogonally bipolar polytopes with poles $\sigma,\nu$ (large blue dots):
the square (left, a $2$-HC) and the cube (right, a $3$-HC) in an oblique
projection. In both figures $\sigma\nu$ (densely dashed) is a diameter of the
circumscribed Thales sphere. Left: the faces $[\sigma,\epsilon_1]$ and
$[\epsilon_1,\nu]$ are segments, so their perpendicularity is the right angle at
$\epsilon_1$, and likewise at $\epsilon_2$. Right: the atoms at $\sigma$ are
$\epsilon_1,\epsilon_2,\epsilon_3$ ($\epsilon_3$ hidden, in grey). Condition (O3)
at the atom $a=\epsilon_1$ says that the edge $[\sigma,a]$ is perpendicular to the
plane of the facet $[a,\nu]$ (shaded blue, the darker of the two shaded facets); this is marked by the right angles
between $[\sigma,a]$ and the two edges of $[a,\nu]$ at $a$. For the coatom $c$, the
upper neighbour of $a$, the facet $[\sigma,c]$ (shaded orange, the lighter one) is perpendicular to
the edge $[c,\nu]$, marked by the right angles between $[c,\nu]$ and the two edges
of $[\sigma,c]$ at $c$ (Remark~\ref{rem:obp-prep}). Right angles between
non-frontal lines appear skewed in the projection.}
\label{fig:obp}
\end{figure}
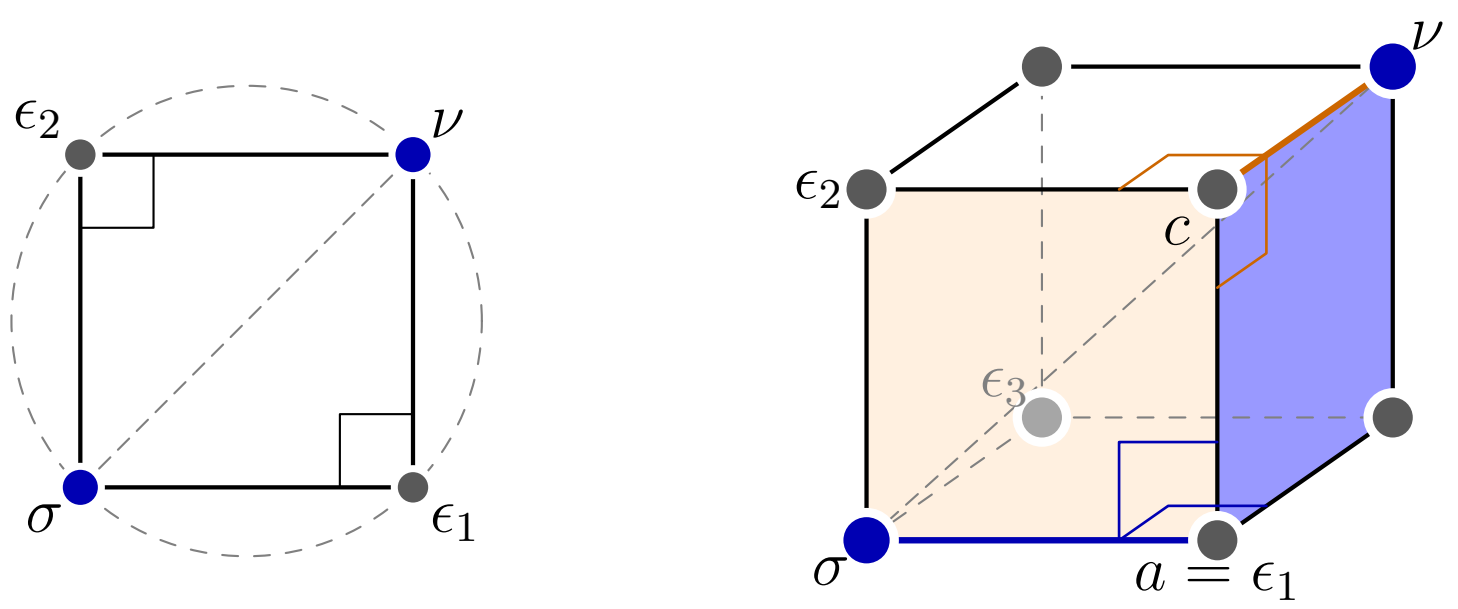

\begin{lm}[full complementarity]\label{lm:full}
Let $W\subseteq\R^n$ be an orthogonally bipolar polytope with poles $\sigma,\nu$. Then
\[
[\sigma,x]\perp[x,\nu]\qquad\text{for every vertex }x\text{ of }W
\text{ with } x\ne\sigma,\nu .
\]
\end{lm}
\begin{proo}
Let $x\ne\nu$ be a vertex of $W$ and let $a_1,\dots,a_r$ be the atoms of
the vertex lattice $(V(W),\le)$ with $a_i\le x$. The edges of $W$ at
$\sigma$ are the faces $[\sigma,a]$ with $a$ an atom; those lying in the
face $[\sigma,x]$ are exactly $[\sigma,a_1],\dots,[\sigma,a_r]$, and their
directions span $U_{[\sigma,x]}$. Hence
\[
U_{[\sigma,x]}=\lin\{a_i-\sigma:\ i=1,\dots,r\}.
\]
Dually, the edges of $[x,\nu]$ at $\nu$ are the faces $[c,\nu]$ with $c$ a
coatom of $(V(W),\le)$, $c\ge x$, so
\[
U_{[x,\nu]}=\lin\{\nu-c:\ c\text{ a coatom},\ c\ge x\}.
\]
Now take $i\in\{1,\dots,r\}$ and a coatom $c\ge x$. Then $a_i\le x\le c$,
hence $[c,\nu]\subseteq[a_i,\nu]$ and $\nu-c\in U_{[a_i,\nu]}$. By (O3) the
faces $[\sigma,a_i]$ and $[a_i,\nu]$ are perpendicular, so
$a_i-\sigma\perp\nu-c$. As the $a_i-\sigma$ span $U_{[\sigma,x]}$ and the
$\nu-c$ span $U_{[x,\nu]}$, the two direction spaces are perpendicular.
\end{proo}

\begin{lm}[heredity]\label{lm:obpface}
Every face of an orthogonally bipolar polytope is again an orthogonally bipolar polytope: if $W$ is an orthogonally bipolar polytope with poles
$\sigma,\nu$ and $x\le y$ in $(V(W),\le)$, then the face $[x,y]$ of $W$ is
an orthogonally bipolar polytope with poles $x,y$.
\end{lm}
\begin{proo}
By (F2), the face $[x,y]$ is a geometric bipolar polytope with poles $x,y$,
whose vertex lattice is the interval $[x,y]$; so (O1) and (O2) hold for it. For (O3) let $a$ be an atom of the
vertex lattice $[x,y]$ of that face, i.e.\ an element of $[x,y]$ covering
$x$. Then
\[
[x,a]\subseteq[\sigma,a],\qquad [a,y]\subseteq[a,\nu]
\]
(since $\sigma\le x\le a\le y\le\nu$), so
$U_{[x,a]}\subseteq U_{[\sigma,a]}$ and
$U_{[a,y]}\subseteq U_{[a,\nu]}$; these two are perpendicular by
Lemma~\ref{lm:full}. Hence $[x,a]\perp[a,y]$.
\end{proo}

\begin{theo}[Thales's theorem for orthogonally bipolar polytopes]
\label{theo:obpthales}
Let $W\subseteq\R^n$ be an  orthogonally bipolar polytope with poles $\sigma,\nu$,
let $T$ be the midpoint of the segment $\sigma\nu$ and
$r=\frac12\delta(\sigma,\nu)$. Then $\delta(T,K)=r$ for every vertex $K$
of $W$; equivalently, the $(n-1)$-dimensional sphere with centre $T$ and
radius $r$ (the sphere with diameter $\sigma\nu$) is circumscribed to $W$.
More generally, for every face $[x,y]$ of $W$ with $x<y$ all vertices
of $[x,y]$ lie on the sphere with diameter $[x,y]$.
\end{theo}
\begin{proo}
For $K\in\{\sigma,\nu\}$ there is nothing to prove; let $K\ne\sigma,\nu$.
The vector $K-\sigma$ is a vertex vector of the face $[\sigma,K]$ and
$\nu-K$ is a vertex vector of $[K,\nu]$, so Lemma~\ref{lm:full} gives
$\langle K-\sigma,\,K-\nu\rangle=0$. Put $a=K-T$ and $b=T-\sigma=\nu-T$,
so $\|b\|=r$; then $K-\sigma=a+b$, $K-\nu=a-b$ and
\[
0=\langle K-\sigma,\,K-\nu\rangle=\|a\|^2-\|b\|^2=\delta(T,K)^2-r^2 .
\]
The statement about a face $[x,y]$ follows by applying this to the orthogonally bipolar polytope
$[x,y]$ with poles $x,y$ (Lemma~\ref{lm:obpface}).
\end{proo}
In other words, the theorem is the vertexwise Thales condition
$\angle\sigma K\nu=\frac\pi2$ supplied by Lemma~\ref{lm:full}; neither an
induction nor the dimension of $W$ is needed.

\begin{coro}\label{cor:obpdiam}
For an orthogonally bipolar polytope $W$ with poles $\sigma,\nu$ the segment $\sigma\nu$ is a
diameter of $W$: $\delta(x,y)\le\delta(\sigma,\nu)$ for all vertices
$x,y$, with equality iff $x$ and $y$ are antipodal on the sphere of
Theorem~\ref{theo:obpthales}, i.e.\ iff $x+y=\sigma+\nu$.
\end{coro}
\begin{proo}
By Theorem~\ref{theo:obpthales} all vertices lie on the sphere with centre
$T$ and radius $r$, so $\delta(x,y)\le2r=\delta(\sigma,\nu)$ by the
triangle inequality; equality holds iff $x-T$ and $T-y$ have the same
direction.
\end{proo}

\begin{rem}[the converse fails]\label{rem:converse}
In the plane the converse of Theorem~\ref{theo:obpthales} holds: a convex
quadrilateral $A,B,C,D$ inscribed in a circle with diameter $AC$ has right
angles at $B$ and $D$, so it is a Hjelmslev quadrilateral. In dimension
three, a geometric bipolar polytope inscribed in the sphere with diameter
$\sigma\nu$ need not be orthogonally bipolar. Let $a=2-\sqrt3$,
$\rho=\sqrt{1-a^2}=\sqrt{4\sqrt3-6}$, $\sigma=(0,0,-1)$, $\nu=(0,0,1)$
and, for $i=0,1,2$,
\[
h_i=\bigl(\rho\cos\tfrac{2\pi i}3,\ \rho\sin\tfrac{2\pi i}3,\ -a\bigr),\qquad
t_i=\bigl(\rho\cos(\tfrac{2\pi i}3+\tfrac\pi2),\
\rho\sin(\tfrac{2\pi i}3+\tfrac\pi2),\ a\bigr).
\]
All eight points lie on the unit sphere, on which $\sigma$ and $\nu$ are
antipodal. For this value of $a$ the quadruples $\sigma h_it_ih_{i+1}$ and
$\nu t_{i-1}h_it_i$ (indices mod $3$) are coplanar, and
$W=\conv\{\sigma,\nu,h_i,t_i\}$ is a combinatorial cube with exactly these
six facets. Each facet contains exactly one pole, so $W$ is a geometric
bipolar polytope with vertex lattice $B_3$, inscribed in the sphere with
diameter $\sigma\nu$. However
$\langle h_0-\sigma,\,t_0-h_0\rangle=2\sqrt3-4\ne0$, so the edge
$[\sigma,h_0]$ is not perpendicular to the facet $[h_0,\nu]=\nu t_2h_0t_0$
and (O3) fails. Rotating the points $t_i$ by another angle $\theta$ instead
of $\frac\pi2$ and adjusting $a$ gives further examples; $\theta=\frac\pi3$
gives the cube.
\end{rem}

\subsection{The Boolean case: Hjelmslev
\texorpdfstring{$n$}{n}-cuboids}\label{sec:boolean}

An \emph{$n$-cuboid} is a polytope combinatorially equivalent to the unit
$n$-cube. Inside the class of orthogonally bipolar polytopes the Boolean
vertex lattices are singled out by the degree of a pole.

\begin{prop}\label{prop:boolean}
Let $W$ be an $n$-dimensional orthogonally bipolar polytope with poles $\sigma,\nu$ and vertex
lattice $(V(W),\le)$. Then the following are equivalent:
\begin{enumerate}[label=(\roman*)]
\item $(V(W),\le)\cong B_n$, the Boolean lattice of rank $n$;
\item $(V(W),\le)$ has exactly $n$ atoms;
\item $\sigma$ is joined by an edge to exactly $n$ vertices, i.e.\ $\sigma$
      is a simple vertex of $W$.
\end{enumerate}
In that case $\F(W)\cong\widetilde I(B_n)$, the lattice of intervals of
$B_n$ with a new least element $\emptyset$, is the face lattice of the
$n$-cube \cite{MetropolisRota}, so $W$ is an $n$-cuboid, and the poles are antipodal. Conversely,
a cubical orthogonally bipolar polytope has antipodal poles by (O1), hence satisfies (i). An orthogonally bipolar polytope with
these properties is called a \emph{Hjelmslev $n$-cuboid} ($n$-HC).
\end{prop}
\begin{proo}

(i)$\Rightarrow$(ii) is clear, as $B_n$ has $n$ atoms.
(ii)$\Rightarrow$(i): by Remark~\ref{rem:obp-prep},
$(V(W),\le)\cong\F(W/\sigma)$, whose atoms are the vertices of the
$(n-1)$-polytope $W/\sigma$. An $(n-1)$-polytope with exactly $n$ vertices
is a simplex, and the face lattice of the $(n-1)$-simplex is $B_n$. (ii)$\Leftrightarrow$(iii): by Remark~\ref{rem:obp-prep}
the atoms of $(V(W),\le)$ are exactly the edges at $\sigma$, so their
number is the degree of $\sigma$. (i)$\Rightarrow$ cube: the faces of $W$
are the intervals of $(V(W),\le)$, so $\F(W)$ is the lattice of intervals
of $B_n$, i.e.\ the face lattice of the $n$-cube. Conversely, if $W$ is an
$n$-cuboid, then by (O1) the poles lie in opposite facets and are
therefore antipodal, and the vertex figure at $\sigma$ is the
$(n-1)$-simplex, so $(V(W),\le)\cong\F(W/\sigma)\cong B_n$.
\end{proo}

\begin{figure}[htbp]
\centering
\begin{tikzpicture}[scale=1.45,line join=round,line cap=round]
  \coordinate (S) at (0.000,0.000);
  \coordinate (C) at (2.553,0.511);
  \coordinate (D) at (0.743,2.512);
  \coordinate (E) at (-0.492,1.832);
  \coordinate (F) at (2.636,2.419);
  \coordinate (G) at (1.648,1.874);
  \coordinate (H) at (0.200,3.475);
  \coordinate (N) at (1.869,3.237);
  \fill[blue!10] (S)--(C)--(G)--(E)--cycle;
  \fill[red!8]   (E)--(G)--(N)--(H)--cycle;
  \draw[black!60,dashed,line width=0.8pt] (S)--(D) (D)--(F) (D)--(H);
  \draw[line width=1.1pt] (S)--(C) (S)--(E) (C)--(F) (C)--(G)
        (E)--(H) (E)--(G) (F)--(N) (G)--(N) (H)--(N);
  \draw[gray,densely dashed] (S)--(N);
  \rightangle{E}{S}{G}
  \node[below left,inner sep=1pt,outer sep=4pt]  at (S) {$\sigma$};
  \node[above right,inner sep=1pt,outer sep=4pt] at (N) {$\nu$};
  \node[below right,inner sep=1pt,outer sep=4pt] at (C) {$C$};
  \node[below left,inner sep=1pt,outer sep=4pt]  at (D) {$D$};
  \node[left,inner sep=1pt,outer sep=4pt]        at (E) {$E$};
  \node[right,inner sep=1pt,outer sep=4pt]       at (F) {$F$};
  \node[below,inner sep=1pt,outer sep=4pt] at (G) {$G$};
  \node[above left,inner sep=1pt,outer sep=4pt]  at (H) {$H$};
  \pole{S} \pole{N}
  \gen{C} \gen{E} \gen{F} \gen{G} \gen{H} \genw{D}
\end{tikzpicture}
\caption{A Hjelmslev $3$-cuboid which is not a rectangular box. In coordinates
$(x_0,y)\in\R\times\R^2$, with $v_1,v_2,v_3$ the vertices of an equilateral triangle
centred at $O$ with $\|v_i\|=1$, its vertices are the poles $\sigma=O$ and
$\nu=(1,O_2)$, the atoms $C,D,E=\frac12(1,v_i)$, $i=1,2,3$, and the points
$F,G,H=\frac45\bigl(1,\frac12(v_i+v_j)\bigr)$ for $\{i,j\}=\{1,2\},\{1,3\},\{2,3\}$, which
are the joins $C\vee D$, $C\vee E$, $D\vee E$ in the vertex lattice; the drawing is an
exact parallel projection. The edges at $\sigma$ make pairwise angles with cosine
$\frac14$, yet (O3) holds; the right angle at $E$ between $\sigma E$ and the facet
$[E,\nu]$ is marked. The shaded facets are the south facet $[\sigma,G]=\sigma CGE$ and
the north facet $[E,\nu]=EG\nu H$; they share the edge $EG$, and each contains exactly
one pole. $D$ is the hidden vertex and its three edges are dashed. A parallelepiped
cannot serve as such an example: if opposite edges are parallel, (O3) forces the
three edges at $\sigma$ to be pairwise orthogonal.}
\label{fig:3hc}
\end{figure}
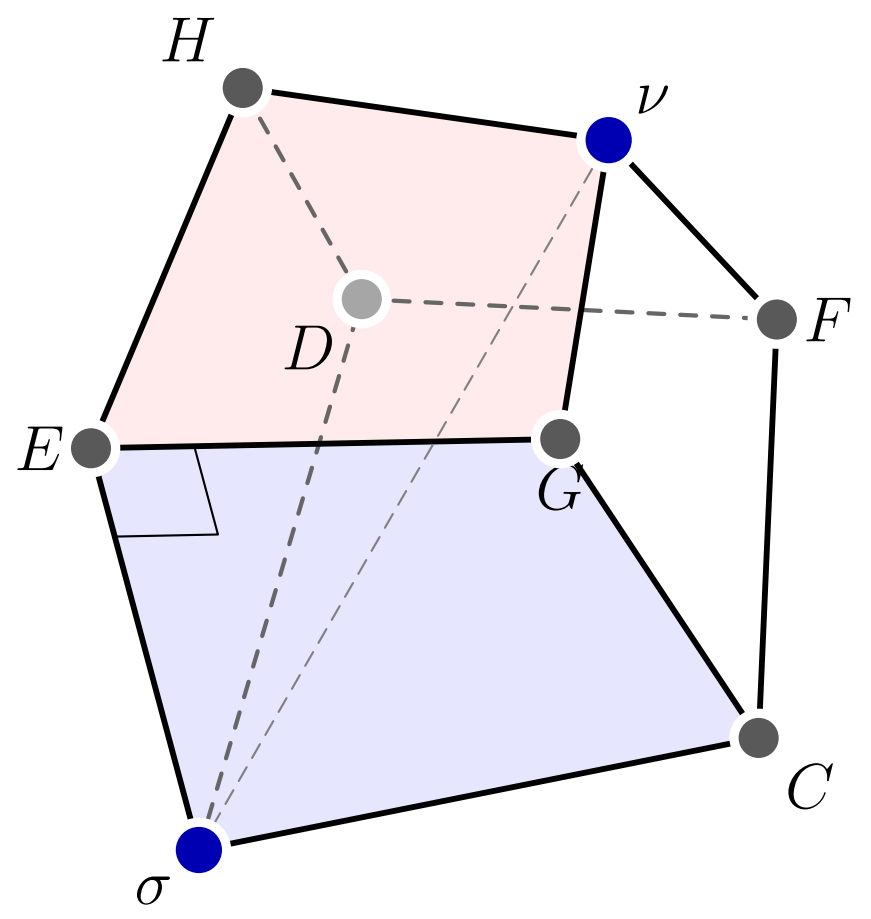

Lemmas~\ref{lm:full}, \ref{lm:obpface} and
Theorem~\ref{theo:obpthales} apply to every $n$-HC. In particular every
$k$-face $[x,y]$ of an $n$-HC, being itself a $k$-HC, has its vertices on
the $(k-1)$-sphere in $\aff[x,y]$ with diameter $[x,y]$; this holds
verbatim for every orthogonally bipolar polytope. A Hjelmslev $3$-cuboid need not be a rectangular box; an
example is shown in Figure~\ref{fig:3hc}.

\section{Orthogonal cutouts}\label{sec:cutouts}

Theorem~\ref{theo:obpthales} starts from an orthogonally bipolar polytope and shows
that its vertices lie on a sphere. In this section we go the other way. The orthogonal
projections of a point onto the faces of a polytope through a vertex always lie on a
Thales sphere (\S\ref{sec:proj}); orthogonally bipolar polytopes are exactly the
cutouts formed by such projections (\S\ref{sec:obp-cutouts}); and this yields a
description of orthogonally bipolar polytopes by a sphere, a pole and lines through
it (\S\ref{sec:lines}).

\subsection{Orthogonal projections}\label{sec:proj}

\begin{lm}\label{lm:orth-thales}
Let $P\subseteq\R^n$ be a polytope, $B\in\R^n$ a point, $V$ a vertex of $P$ and, for a face $G$ of $P$ with $V\in G$, let
$x_G$ be the orthogonal projection of $B$ onto $\aff G$; put
$x_{\{V\}}=V$ and $x_P=B$. Then
\[
\langle x_G-V,\,x_G-B\rangle=0
\]
for every face $G\ni V$; equivalently, all the points $x_G$ lie on the
sphere with diameter $[V,B]$.
\end{lm}
\begin{proo}
$x_G-B\perp U_G$, while $x_G-V\in U_G$; hence the two vectors are
orthogonal, which is the $n$-dimensional Thales condition for the sphere
with diameter $[V,B]$.
\end{proo}

\begin{figure}[htbp]
\centering
\begin{tikzpicture}[scale=0.95,line cap=round,line join=round]
  \coordinate (A)  at (0,0);
  \coordinate (C)  at (4.5,0);
  \coordinate (E)  at (1,3.4);
  \coordinate (B)  at (1.5,1.1);
  \coordinate (y1) at (1.5,0);
  \coordinate (y2) at (0.4172,1.4185);
  \fill[gray!8] (A)--(C)--(E)--cycle;
  \fill[blue!30,fill opacity=0.7] (A)--(y1)--(B)--(y2)--cycle;
  \draw[thick] (A)--(C)--(E)--cycle;
  \draw[thick] (A)--(y1)--(B)--(y2)--cycle;
  \draw[gray!75,dotted] (B)--(y1) (B)--(y2);
  \draw[gray!85,dashed,line width=0.7pt] (0.75,0.55) circle (0.9301);
  \draw[gray!60,dotted] (A)--(B);
  \rightangle{y1}{A}{B}
  \rightangle{y2}{A}{B}
  \pole{A} \pole{B} \gen{y1} \gen{y2} \gen{C} \gen{E}
  \node[below left, fill=white,inner sep=1pt] at (A)  {$A$};
  \node[right,      fill=white,inner sep=1pt] at (B)  {$B$};
  \node[below,      fill=white,inner sep=1pt] at (y1) {$y_1$};
  \node[left,       fill=white,inner sep=1pt] at (y2) {$y_2$};
  \node[below,      fill=white,inner sep=1pt] at (C)  {$C$};
  \node[left,       fill=white,inner sep=1pt] at (E)  {$E$};
\end{tikzpicture}\hspace{88.6pt}%
\begin{tikzpicture}[scale=0.78,line cap=round,line join=round]
  \coordinate (A)  at (0,0);
  \coordinate (C)  at (4,0);
  \coordinate (E)  at (-3,4);
  \coordinate (B)  at (-1.2,2.5);
  \coordinate (x1) at (-1.2,0);
  \coordinate (x2) at (-1.632,2.176);
  \fill[gray!8] (A)--(C)--(E)--cycle;
  \draw[thick] (A)--(C)--(E)--cycle;
  \draw[gray!75,dotted] (B)--(x1) (B)--(x2);
  \draw[gray!50,densely dotted] (x1)--(A);
  \draw[gray!85,dashed,line width=0.7pt] (-0.6,1.25) circle (1.3865);
  \draw[gray!60,dotted] (A)--(B);
  \rightangle{x1}{A}{B}
  \rightangle{x2}{A}{B}
  \pole{A} \pole{B} \gen{x1} \gen{x2} \gen{C} \gen{E}
  \node[below left, fill=white,inner sep=1pt] at (A)  {$A$};
  \node[right,      fill=white,inner sep=1pt] at (B)  {$B$};
  \node[below,      fill=white,inner sep=1pt] at (x1) {$x_1$};
  \node[below left, fill=white,inner sep=1pt] at (x2) {$x_2$};
  \node[below right,fill=white,inner sep=1pt] at (C)  {$C$};
  \node[below left, fill=white,inner sep=1pt] at (E)  {$E$};
\end{tikzpicture}
\caption{Left: the orthogonal configuration of $B$ at the vertex $A$ of
the acute triangle $ACE$. Both orthogonal projections of $B$ onto the edges
through $A$ fall into their relative interiors, so the cutout
$A,y_1,B,y_2$ (shaded) is admissible: it has right angles at the
generators $y_1,y_2$ and the circle with diameter $AB$ (dashed) is
circumscribed to it. Right: for the obtuse triangle $ACE$ the orthogonal
projections of $B$ at $A$ do not form a configuration. The perpendicular from $B$ to the line $AC$ meets
it at $x_1$, at distance $1.2$ \emph{outside} the edge $AC$ (the missed
part of the base is drawn dotted), while the projection $x_2$ onto $AE$ is
interior. Hence $x_1\notin\relint AC$ and the cutout at $A$ is not
admissible; both $x_1$ and $x_2$ nevertheless lie on the sphere with
diameter $AB$: the Thales sphere always exists, it is the relative
interiors that may fail.}
\label{fig:cover}
\end{figure}
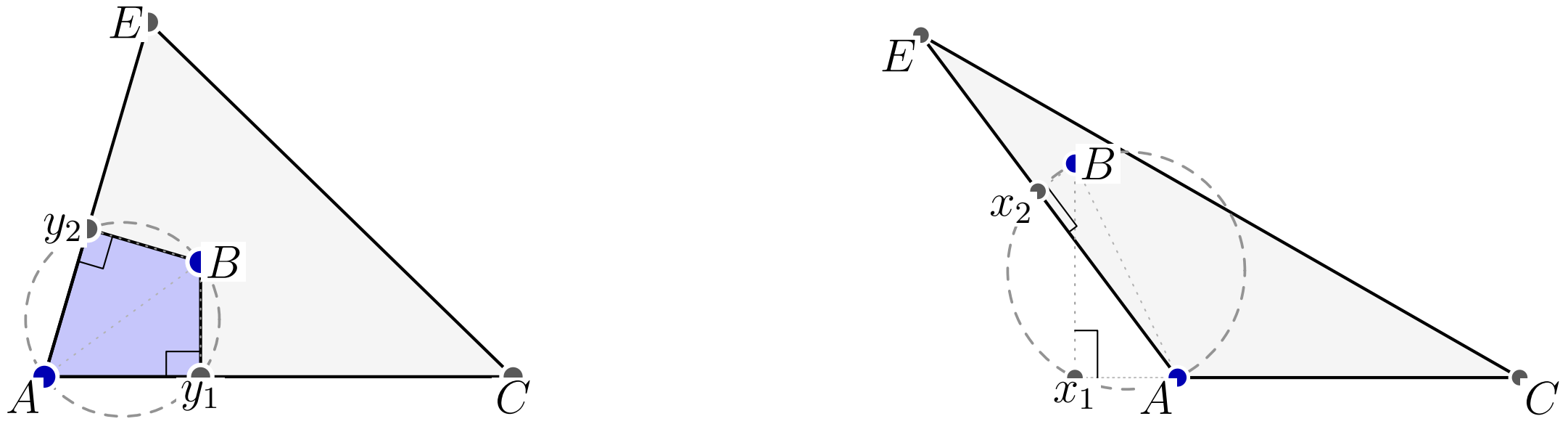

The projections need not lie in the relative interiors of their faces; in that case
they do not form a configuration (Figure~\ref{fig:cover}).

\subsection{Orthogonally bipolar polytopes are orthogonal cutouts}\label{sec:obp-cutouts}

\begin{prop}\label{prop:cut-obp}
Let $P\subseteq\R^n$ be an $n$-polytope, $V$ a vertex of $P$ and $B\in\inte P$, and
suppose that the orthogonal projections $x_G$ of $B$ onto the affine hulls of the
faces $G\ni V$ form the orthogonal configuration of $B$ at $V$ and that its cutout
$C=\conv\{x_G:\ V\in G\}$ is admissible (\S\ref{sec:notation}). Then $C$ is an
orthogonally bipolar polytope with poles $V$ and $B$.
\end{prop}
\begin{proo}
$C$ is a geometric bipolar polytope with poles $V,B$ by admissibility, so it remains
to check (O3). For an edge $e$ of $P$ at $V$ and a facet $\Phi\supseteq e$ we have
$x_e-V\in U_e\subseteq U_\Phi$ and $x_\Phi-B\perp U_\Phi$. The vectors $B-x_\Phi$,
$\Phi\supseteq e$ a facet of $P$, are the edges of the face $[x_e,B]$ of $C$ at $B$
and so span its direction space; hence the faces $[V,x_e]$ and $[x_e,B]$ of $C$ are
perpendicular.
\end{proo}

The converse says that every orthogonally bipolar polytope arises in this way.

\begin{prop}\label{prop:obp-cut}
Let $W$ be an orthogonally bipolar polytope with poles $\sigma,\nu$. Then every vertex
$x$ of $W$ is the orthogonal projection of $\nu$ onto $\aff[\sigma,x]$, and it lies in
the relative interior of the face of the tangent cone of $W$ at $\sigma$ spanned by
$[\sigma,x]$. Consequently, if $P$ is the tangent cone of $W$ at $\sigma$ cut off by
a hyperplane perpendicular to $\sigma\nu$ beyond $\nu$, then $W$ is the cutout at
$\sigma$ of the orthogonal configuration of $\nu$ in $P$.
\end{prop}
\begin{proo}
By Lemma~\ref{lm:full}, $[\sigma,x]\perp[x,\nu]$, so the vector $\nu-x$ is orthogonal
to the direction space of $\aff[\sigma,x]$; as $x\in\aff[\sigma,x]$, the point $x$ is
the orthogonal projection of $\nu$ onto $\aff[\sigma,x]$. The face $[\sigma,x]$ is the
smallest face of $W$ containing $\sigma$ and $x$, so $x$ lies in no smaller face of the
tangent cone. Every atom $a$ lies on the sphere with diameter $[\sigma,\nu]$
(Theorem~\ref{theo:obpthales}), i.e.\ $\langle a-\sigma,a-\nu\rangle=0$ (\S\ref{sec:notation});
writing $a-\nu=(a-\sigma)-(\nu-\sigma)$, this gives
$\langle a-\sigma,\nu-\sigma\rangle=\|a-\sigma\|^2>0$, and $\nu$ lies in no facet of $W$ containing $\sigma$;
so $\nu$ is an interior point of the tangent cone and of $P$, and the faces of $P$
containing $\sigma$ have the same affine hulls as the corresponding faces of the
tangent cone.
\end{proo}

\begin{rem}\label{rem:thm-vs-lemma}
By Proposition~\ref{prop:obp-cut}, Theorem~\ref{theo:obpthales} is the case $P=W$,
$V=\sigma$, $B=\nu$ of Lemma~\ref{lm:orth-thales}: the vertices of an orthogonally
bipolar polytope are the orthogonal projections of $\nu$ onto the affine hulls of
the faces containing $\sigma$, and Lemma~\ref{lm:orth-thales} puts them on the sphere
with diameter $[\sigma,\nu]$. Lemma~\ref{lm:orth-thales} is more general: it holds for
every polytope, every vertex and every point, although in general the projections
need not even form a configuration (Figure~\ref{fig:cover}, right), let alone the
vertex set of an orthogonally bipolar polytope.
\end{rem}

\subsection{A description by a sphere and lines}\label{sec:lines}

Propositions~\ref{prop:cut-obp} and~\ref{prop:obp-cut} show that an orthogonally
bipolar polytope is determined by very little data: its two poles and the lines
through $\sigma$ carrying its edges.

\begin{coro}\label{cor:lines}
Let $\sigma\ne\nu$, let $\Sigma$ be the sphere with diameter $[\sigma,\nu]$, and let
$\ell_1,\dots,\ell_m$ be lines through $\sigma$ meeting $\Sigma$ again in points
$a_1,\dots,a_m$. Let $T=\sigma+\cone\{a_i-\sigma\}$, and assume that $\nu$ is an
interior point of $T$ and that the $a_i$ lie on distinct edges of $T$.
\begin{enumerate}[label=(\alph*)]
\item There is at most one orthogonally bipolar polytope with poles $\sigma,\nu$ whose
      atoms are $a_1,\dots,a_m$, namely the cutout at $\sigma$ of the orthogonal
      configuration of $\nu$ in $T$ (cut off beyond $\nu$).
\item It exists if and only if the section $Q'$ of $T$ by the hyperplane through $\nu$
      perpendicular to $\sigma\nu$ satisfies the interval condition of
      \cite[Definition~5.1]{Bipolar} with respect to the point $\nu$.
\item For $n=3$ this means: for any two lines $\ell_i,\ell_j$ spanning a face of $T$,
      the orthogonal projection of $\nu$ onto their plane lies strictly inside the
      angle between them.
\end{enumerate}
\end{coro}
\begin{proo}
(a) is Proposition~\ref{prop:obp-cut}. (b) is Theorem~\ref{prop:obpreal},
applied with the base $Q'$ and $B=\nu$; the atoms are automatically the projections of
$\nu$ onto the lines $\ell_i$, by Thales's theorem. (c) For a polygon $Q'$ the
interval condition holds iff the foot of $\nu$ on the line of every edge lies in the
relative interior of that edge \cite[Theorem~5.4(a)]{Bipolar} and Proposition~\ref{prop:cutout}(c);
the projection of $\nu$ onto the plane spanned by $\ell_i,\ell_j$ lies on the ray from
$\sigma$ through that foot.
\end{proo}

\begin{figure}[htbp]
\centering
\begin{tikzpicture}[scale=2.0,line cap=round,line join=round]
\draw[gray!60] (0.000,0.767) circle (1);
\fill[red!10] (0.000,0.000)--(0.074,0.825)--(1.560,0.444)--cycle;
\draw[gray!70,thin] (0.074,0.825)--(1.560,0.444)--(-0.446,2.063)--cycle;
\draw[gray!70,densely dotted] (0.074,0.825)--(-0.371,0.939);
\draw[gray!60,dashed] (0.000,0.000)--(0.000,1.535);
\draw[black,line width=0.8pt] (0.000,0.000)--(0.074,0.825);
\draw[black,line width=0.8pt] (0.000,0.000)--(1.560,0.444);
\draw[black,line width=0.8pt] (0.000,0.000)--(-0.446,2.063);
\draw[red!70!black,densely dotted,line width=0.7pt] (0.000,0.000)--(-0.371,0.939);
\draw[red!70!black,dashed,line width=0.6pt] (0.000,1.535)--(-0.297,0.751);
\fill[black] (0.057,0.631) circle (1.5pt);
\fill[black] (0.669,0.190) circle (1.5pt);
\fill[black] (-0.365,1.692) circle (1.5pt);
\fill[red!70!black] (-0.297,0.751) circle (1.7pt);
\fill[gray] (-0.371,0.939) circle (1.2pt);
\fill[blue!70!black] (0.000,0.000) circle (1.9pt); \fill[blue!70!black] (0.000,1.535) circle (1.9pt);
\node[black,font=\small,below=3pt,inner sep=0pt] at (0.000,0.000) {$\sigma$};
\node[black,font=\small,right=3pt,inner sep=0pt] at (0.000,1.535) {$\nu$};
\node[black,font=\small,below right=3pt,inner sep=0pt] at (0.057,0.631) {$a_1$};
\node[black,font=\small,below right=3pt,inner sep=0pt] at (0.669,0.190) {$a_2$};
\node[black,font=\small,left=3pt,inner sep=0pt] at (-0.365,1.692) {$a_3$};
\node[red!70!black,font=\small,left=3pt,inner sep=0pt] at (-0.297,0.751) {$x$};
\node[gray,font=\small,left=3pt,inner sep=0pt] at (-0.371,0.939) {$q$};
\node[black,font=\small,above=3pt,inner sep=0pt] at (0.074,0.825) {$\ell_1$};
\node[black,font=\small,right=3pt,inner sep=0pt] at (1.560,0.444) {$\ell_2$};
\node[black,font=\small,left=3pt,inner sep=0pt] at (-0.446,2.063) {$\ell_3$};
\end{tikzpicture}
\caption{Lines through a pole need not suffice (Corollary~\ref{cor:lines}). Three lines
$\ell_1,\ell_2,\ell_3$ through $\sigma$ meet the sphere with diameter $[\sigma,\nu]$
again in $a_1,a_2,a_3$; the thin triangle is the section of the cone $T$ they span by
the plane through $\nu$ perpendicular to $\sigma\nu$. The orthogonal projection $x$ of
$\nu$ onto the plane of $\ell_1,\ell_2$ lies on the sphere (Lemma~\ref{lm:orth-thales})
and on the ray from $\sigma$ through the foot $q$ of $\nu$ on the line of the edge of
the section between $\ell_1$ and $\ell_2$ (dotted); but $q$ lies outside that edge, so
$x$ lies outside the angle between $\ell_1$ and $\ell_2$ (shaded). Hence no orthogonally
bipolar polytope has the poles $\sigma,\nu$ and the atoms $a_1,a_2,a_3$; the section
violates the interval condition. Coordinates: $\sigma=O$, $\nu=(0,0,2)$, and the lines
pass through $(0.48,-1,2)$, $(2.08,-1,2)$, $(-0.72,0.6,2)$; then
$x=(0,-0.8,1.6)$ and $q=(0,-1,2)$.}
\label{fig:lines}
\end{figure}
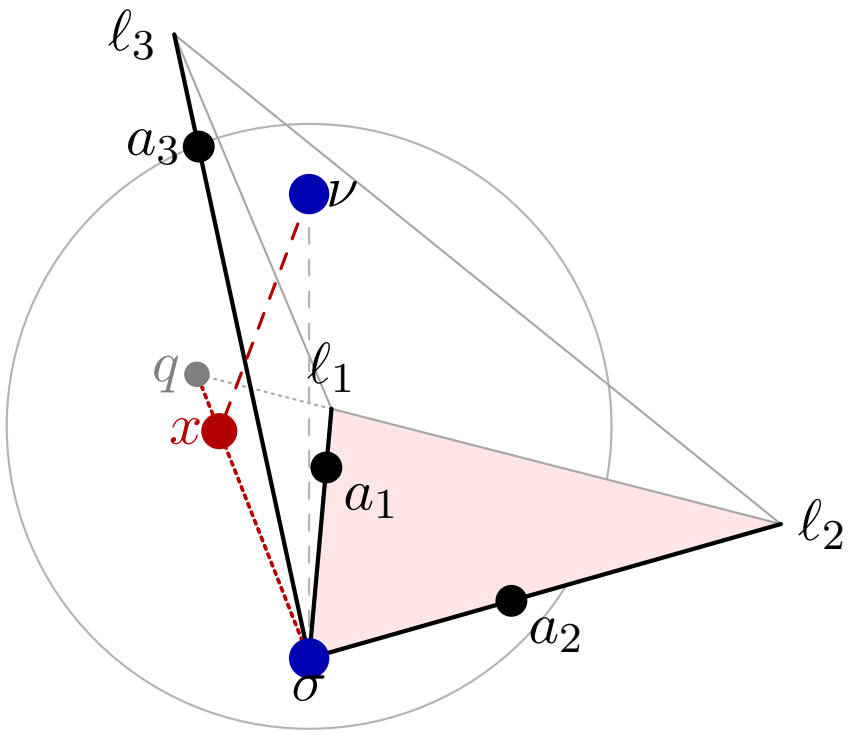

Condition (b) does not hold automatically, even in $\R^3$: Figure~\ref{fig:lines} shows
three lines through $\sigma$ for which the projection of $\nu$ onto the plane of two of
them misses the angle between them.

Thus the essential condition is (b), and it depends on the position of the axis
$\sigma\nu$ inside the cone $T$. For a fixed section $Q'$, the set of points
$O\in\inte Q'$ about which the interval condition holds is an explicit open convex
polytope, the set of points about which $Q'$ is perfectly centered
Theorem~\ref{conj:I}, Proposition~\ref{rem:centres} and Figure~\ref{fig:centres}; in
the plane this is the intersection of the open strips perpendicular to the edges.

\subsection{The apex cutout and orthogonally bipolar realizations}\label{sec:cutout}

In this subsection $Q\subset\R^d$ is a polytope with $O\in\inte Q$, $Q^\circ$ is its polar,
and $W(\F(Q))$ is the bipolar polytope of \cite{Bipolar} whose vertices are the faces of
$Q$ and whose faces are the intervals of $\F(Q)$, i.e.\ the interval polytope of $Q$. Its
convex realization as the polar of the Cayley polytope $\mathrm{Cay}(Q,Q^\circ)$, under
the interval condition (I) of \cite[Definition~5.1]{Bipolar}, is a theorem of Broadie
\cite{Broadie} \cite[Theorem~5.2]{Bipolar}. Here we show that this realization has a direct
metric meaning: it can be chosen orthogonally bipolar, as the cutout of the apex of a
pyramid by orthogonal projections. For a
nonempty face $F$ of $Q$ let $q_F$ be the point of $\aff F$ nearest to $O$ (so
$q_{\{u\}}=u$).

\begin{de}[orthogonal reciprocity]\label{def:star}
$Q$ is \emph{orthogonally reciprocal} if for every nonempty proper face $F$:
\begin{enumerate}[label=($\star$\arabic*)]
\item $q_F\in\relint F$;
\item $\langle q_F,u\rangle<\|q_F\|^2$ for every vertex $u\notin F$.
\end{enumerate}
\end{de}

Since $\langle q_F,u\rangle=\|q_F\|^2$ for $u\in\aff F$, the conditions say that
$q_F\in\relint F\cap\relint N_Q(F)$: the direction of the foot of the perpendicular
from $O$ to $\aff F$ points into $F$ and exposes $F$. They are invariant under
polarity, and regular polytopes satisfy them.

\begin{lm}\label{lem:half}
If $Q$ is orthogonally reciprocal, then the implication ``$\Leftarrow$'' of (I)
holds.
\end{lm}
\begin{proo}
Let $\emptyset\ne H\subseteq G\ne Q$ and put $c=q_H+\varepsilon q_G$ with $\varepsilon>0$.
By ($\star$), $q_H\in\relint N_Q(H)$, and $q_G\in N_Q(G)\subseteq N_Q(H)$;
so $c\in\relint N_Q(H)$. Also $q_G\in\relint G$ gives
$q_G\in\relint\cone(G)$, and $q_H\in H\subseteq G$ gives
$q_H\in\cone(G)$; so $c\in\relint\cone(G)$.
\end{proo}

From now on until the end of this subsection, consider the pyramid over $Q$. Fix
$\varepsilon,\beta>0$, identify $\R^{d+1}=\R\times\R^d$, let $T=\cone\{a_u:u\in V(Q)\}$ with
$a_u=(1,\varepsilon u)$, and let $B=(\beta,O_d)$. For a nonempty proper face $F$ of $Q$ let
$x_F$ be the orthogonal projection of $B$ onto the linear span of the face
$T_F=\cone\{a_u:u\in F\}$ of $T$, and put $x_\emptyset=O$, $x_Q=B$. Let
$H_u=\{x:\langle x-B,a_u\rangle=0\}$ be the hyperplane through $B$ orthogonal to the
edge of $T$ through $a_u$, and
\[
X=T\cap\bigcap_{u\in V(Q)}\{x:\langle x-B,a_u\rangle\le0\}.
\]
If the cone is cut by a hyperplane $x_0=h$ with $h>\beta$, then $T\cap\{x_0\le h\}$ is a
pyramid over a homothetic copy of $Q$, $B$ is an interior point of it, and the $x_F$
are the orthogonal projections of $B$ onto the faces at the apex; $\conv\{x_F\}$ is
the \emph{apex cutout}. In the terms of \S\ref{sec:notation}: the $x_F$ are the
orthogonal projections of $B$ onto the affine hulls of the faces of the pyramid that
contain the apex (the face over $F$ has the same affine hull as $T_F$); if they lie in
the relative interiors of these faces, they form the orthogonal configuration of $B$ at
the apex, and $\conv\{x_F\}$ is its cutout.

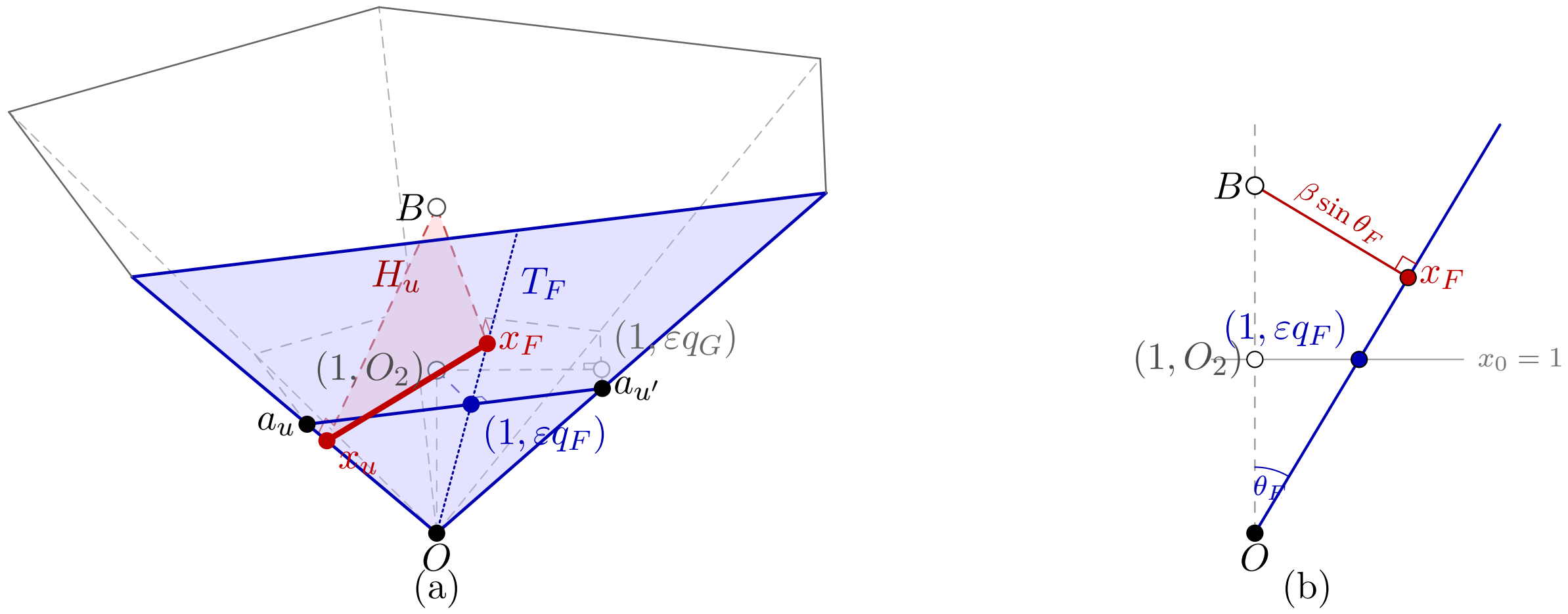
\begin{figure}[htbp]
\centering
\begin{tikzpicture}[scale=1.75,line join=round,line cap=round]
\draw[black!45,densely dashed,line width=0.4pt] (0.000,0.000)--(2.207,2.731);
\draw[black!45,densely dashed,line width=0.4pt] (0.000,0.000)--(-0.332,3.027);
\draw[black!45,densely dashed,line width=0.4pt] (0.000,0.000)--(-2.462,2.424);
\draw[black!50,dashed,line width=0.45pt] (0.953,0.833)--(0.939,1.162);
\draw[black!50,dashed,line width=0.45pt] (0.939,1.162)--(-0.141,1.288);
\draw[black!50,dashed,line width=0.45pt] (-0.141,1.288)--(-1.048,1.031);
\draw[black!50,dashed,line width=0.45pt] (-1.048,1.031)--(-0.746,0.627);
\draw[black!45,dashed] (0.000,0.000)--(0.000,1.876);
\draw[blue!60!black,dashed,line width=0.5pt] (0.000,0.938)--(0.198,0.742);
\draw[black!45,dashed,line width=0.5pt] (0.000,0.938)--(0.949,0.943);
\draw[blue!60!black,line width=0.4pt] (0.165,0.775)--(0.260,0.786)--(0.293,0.753);
\draw[black!50,line width=0.4pt] (0.849,0.942)--(0.847,0.977)--(0.947,0.978);
\fill[red!14] (0.000,1.876)--(-0.632,0.532)--(0.291,1.091)--cycle;
\draw[red!70!black,dashed,line width=0.6pt] (0.000,1.876)--(0.291,1.091) (0.000,1.876)--(-0.632,0.532);
\draw[red!70!black,line width=0.4pt] (0.263,1.167)--(0.280,1.231)--(0.308,1.155);
\draw[red!70!black,line width=0.4pt] (-0.590,0.620)--(-0.639,0.661)--(-0.681,0.573);
\filldraw[fill=white,draw=black,line width=0.5pt] (0.000,1.876) circle (1.4pt);
\filldraw[fill=white,draw=black,line width=0.5pt] (0.000,0.938) circle (1.4pt);
\filldraw[fill=white,draw=black!55,line width=0.5pt] (0.949,0.943) circle (1.3pt);
\fill[blue!20,opacity=0.55] (0.000,0.000)--(-1.753,1.475)--(2.240,1.958)--cycle;
\fill[white,opacity=0.3] (-1.753,1.475)--(2.240,1.958)--(2.207,2.731)--(-0.332,3.027)--(-2.462,2.424)--cycle;
\draw[blue!70!black,line width=0.9pt] (-0.746,0.627)--(0.953,0.833);
\draw[blue!60!black,densely dotted,line width=0.6pt] (0.000,0.000)--(0.466,1.743);
\draw[red!75!black,line width=1.7pt] (0.291,1.091)--(-0.632,0.532);
\draw[blue!70!black,line width=0.9pt] (0.000,0.000)--(-1.753,1.475);
\draw[blue!70!black,line width=0.9pt] (0.000,0.000)--(2.240,1.958);
\draw[blue!70!black,line width=0.9pt] (-1.753,1.475)--(2.240,1.958);
\draw[black!60,line width=0.5pt] (2.240,1.958)--(2.207,2.731);
\draw[black!60,line width=0.5pt] (2.207,2.731)--(-0.332,3.027);
\draw[black!60,line width=0.5pt] (-0.332,3.027)--(-2.462,2.424);
\draw[black!60,line width=0.5pt] (-2.462,2.424)--(-1.753,1.475);
\fill[black] (0.000,0.000) circle (1.4pt);
\fill[black] (-0.746,0.627) circle (1.4pt);
\fill[black] (0.953,0.833) circle (1.4pt);
\fill[blue!70!black] (0.198,0.742) circle (1.4pt);
\fill[red!75!black] (0.291,1.091) circle (1.4pt);
\fill[red!75!black] (-0.632,0.532) circle (1.4pt);
\node[black,font=\small,inner sep=3pt,below] at (0.000,0.000) {$O$};
\node[black,font=\small,inner sep=3pt,left] at (0.000,1.876) {$B$};
\node[red!75!black,font=\small,inner sep=3pt,right] at (0.291,1.091) {$x_F$};
\node[red!75!black,font=\small,inner sep=3pt,below right] at (-0.632,0.532) {$x_u$};
\node[blue!70!black,font=\small,inner sep=3pt,below right] at (0.198,0.742) {$(1,\varepsilon q_F)$};
\node[black!60,font=\small,inner sep=3pt,above right] at (0.949,0.943) {$(1,\varepsilon q_{G})$};
\node[black!70,font=\small,inner sep=3pt,left] at (0.000,0.938) {$(1,O_2)$};
\node[black,font=\small,inner sep=3pt,left] at (-0.746,0.627) {$a_u$};
\node[black,font=\small,inner sep=3pt,right] at (0.953,0.833) {$a_{u'}$};
\node[blue!70!black,font=\small,inner sep=3pt,below] at (0.619,1.593) {$T_F$};
\node[red!60!black,font=\small,inner sep=3pt,left] at (-0.022,1.478) {$H_u$};
\node[font=\small] at (0.000,-0.320) {(a)};
\end{tikzpicture}\hfill
\begin{tikzpicture}[scale=1.75,line join=round,line cap=round]
\draw[black!45,dashed] (0.000,0.000)--(0.000,2.350);
\draw[black!40] (-0.25,1)--(1.200,1);
\node[font=\scriptsize,black!55,anchor=west] at (1.200,1) {$x_0=1$};
\draw[blue!70!black,line width=0.8pt] (0.000,0.000)--(1.410,2.350);
\draw[red!70!black,line width=0.7pt] (0.000,2.000)--(0.882,1.471);
\draw[red!70!black,line width=0.4pt] (0.805,1.517)--(0.851,1.594)--(0.929,1.548);
\draw[blue!70!black,line width=0.4pt] (0,0.38) arc (90:59.04:0.38);
\node[blue!70!black,font=\scriptsize] at (0.092,0.260) {$\theta_F$};
\filldraw[fill=black,draw=black,line width=0.4pt] (0.000,0.000) circle (1.3pt);
\filldraw[fill=white,draw=black,line width=0.4pt] (0.000,1.000) circle (1.3pt);
\filldraw[fill=blue!70!black,draw=black,line width=0.4pt] (0.600,1.000) circle (1.3pt);
\filldraw[fill=red!75!black,draw=black,line width=0.4pt] (0.882,1.471) circle (1.3pt);
\filldraw[fill=white,draw=black,line width=0.5pt] (0.000,2.000) circle (1.4pt);
\node[black,font=\small,inner sep=3pt,below] at (0.000,0.000) {$O$};
\node[black,font=\small,inner sep=3pt,left] at (0.000,2.000) {$B$};
\node[black!70,font=\small,inner sep=3pt,left] at (0.000,1.000) {$(1,O_2)$};
\node[blue!70!black,font=\small,inner sep=3pt,above left] at (0.600,1.000) {$(1,\varepsilon q_F)$};
\node[red!75!black,font=\small,inner sep=3pt,right] at (0.882,1.471) {$x_F$};
\node[red!70!black,font=\scriptsize,anchor=south,rotate=-31.0,inner sep=2pt] at (0.441,1.735) {$\beta\sin\theta_F$};
\node[font=\small] at (0.3,-0.32) {(b)};
\end{tikzpicture}
\caption{Lemma~\ref{lem:tall} for $d=2$, with $\varepsilon=1$. (a) The cone $T$, seen with
the face $T_F$ (shaded blue, translucent) in front; hidden lines are dashed. The cone is
spanned by the points $a_u=(1,\varepsilon u)$, which form the section $\{1\}\times\varepsilon Q$
of the pentagon $Q$ at height $1$, with centre $(1,O_2)$ on the axis. For the edge
$F=[u,u']$ of $Q$, the point $(1,\varepsilon q_F)$ is the foot of the perpendicular from
$(1,O_2)$ to the line of $F$; likewise $(1,\varepsilon q_G)$ for another edge $G$. The
point $x_F$ is the foot of the perpendicular from the interior point $B=(\beta,O_2)$
(open dot) to the plane of $T_F$; it lies on the ray through $(1,\varepsilon q_F)$ (dotted).
The point $x_u$ is the foot of the perpendicular from $B$ to the edge of $T$ through
$a_u$. The triangle $B\,x_u\,x_F$ lies in the plane $H_u$ through $B$ orthogonal to
$a_u$, inside $T$, and meets $T_F$ in the segment $[x_u,x_F]$ (thick).
(b) The section through the axis and $(1,\varepsilon q_F)$. Since $B$ lies on the axis,
neither $B$ nor the normal of the plane of $T_F$ has a component along $F$; so the
perpendicular from $B$ to that plane lies in this section, and its foot $x_F$ lies on
the line through $(1,\varepsilon q_F)$. This line makes the angle $\theta_F$ with the axis,
$\tan\theta_F=\varepsilon\|q_F\|$, and $\|B-x_F\|=\beta\sin\theta_F$.}
\label{fig:tall}
\end{figure}

Figure~\ref{fig:tall} shows the situation for $d=2$.

\begin{lm}\label{lem:tall}
Let $F$ be a nonempty proper face of $Q$ and $m=\|q_F\|^2$. Then
$x_F=s_F(1,\varepsilon q_F)$ with $s_F=\beta/(1+\varepsilon^2m)$, and for every vertex $u$ of $Q$
\[
\langle x_F-B,\,a_u\rangle=\frac{\varepsilon^2\beta}{1+\varepsilon^2m}\bigl(\langle q_F,u\rangle-m\bigr).
\]
In particular:
\begin{enumerate}[label=(\alph*)]
\item $x_F\in H_u$ for every vertex $u\in F$;
\item $x_F\in\relint T_F$ if and only if $q_F\in\relint F$.
\end{enumerate}
\end{lm}
\begin{proo}
Let $U_F$ be the direction space of $\aff F$, so that the points of $\aff F$ are
$q_F+w$ with $w\in U_F$. As usual (\cite[\S2.2]{Bipolar}) we identify a point with
its position vector; thus $q_F$ also denotes the vector from $O$ to the point $q_F$.
Since $q_F$ is the point of $\aff F$ nearest to $O$, this vector is orthogonal to
every $w\in U_F$: $\langle q_F,w\rangle=0$. The linear span of $T_F$ consists of the
vectors $(s,\ \varepsilon(sq_F+w))$ with $s\in\R$ and $w\in U_F$; in particular it contains
the vectors $(0,w)$, $w\in U_F$. Hence the vector $B-x_F$, which is orthogonal to the span, is orthogonal to all
$(0,w)$, $w\in U_F$; so is the vector $B=(\beta,O_d)$, which lies on the axis. Therefore
$x_F=B-(B-x_F)$ is orthogonal to all $(0,w)$, $w\in U_F$. For a vector
$(s,\varepsilon(sq_F+w))$ of the span, this means $\langle sq_F+w,w'\rangle=\langle w,w'\rangle=0$
for all $w'\in U_F$, i.e.\ $w=O_d$. So $x_F=s_F(1,\varepsilon q_F)$ for some $s_F\in\R$: the
perpendicular from $B$ to the span lies in the plane through the axis and
$(1,\varepsilon q_F)$ (Figure~\ref{fig:tall}(b)). Since $(1,\varepsilon q_F)$ lies in the span,
$\langle x_F-B,(1,\varepsilon q_F)\rangle=0$, i.e.\ $(s_F-\beta)+\varepsilon^2s_Fm=0$, which gives
$s_F=\beta/(1+\varepsilon^2m)$; in particular $s_F>0$. Now
$\langle x_F-B,a_u\rangle=(s_F-\beta)+\varepsilon^2s_F\langle q_F,u\rangle
=\varepsilon^2s_F\bigl(\langle q_F,u\rangle-m\bigr)$, which is the stated expression.

(a) For a vertex $u\in F$ we have $u\in\aff F$, and $q_F$ is the point of $\aff F$
nearest to $O$, so $q_F\perp u-q_F$, i.e.\ $\langle q_F,u\rangle=\|q_F\|^2=m$; hence
$\langle x_F-B,a_u\rangle=0$.

(b) The face $T_F=\cone\{a_u:u\in F\}$ consists of the points $s(1,\varepsilon q)$ with
$s\ge0$ and $q\in F$, and its relative interior of those with $s>0$ and
$q\in\relint F$. Since $x_F=s_F(1,\varepsilon q_F)$ with $s_F>0$, the claim follows.
\end{proo}

\begin{prop}[the apex cutout]\label{prop:cutout}
\begin{enumerate}[label=(\alph*)]
\item $X$ is a polytope, and for $0<z_0<\beta$ its polar with respect to $(z_0,O_d)$
      is a Cayley polytope of positive multiples of $Q$ and $Q^\circ$ placed in two
      parallel hyperplanes. In particular
      $\F(X)\cong\F(\mathrm{Cay}(Q,Q^\circ))^{\mathrm{op}}$, the facets of $X$ being
      $X\cap H_u$, $u\in V(Q)$, and $X\cap\aff T_\Phi$, $\Phi$ a facet of $Q$.
\item If $Q$ satisfies (I), then the vertices of $X$ are exactly the pairwise distinct
      points $x_F$, $F\in\F(Q)$, so $X$ is the apex cutout; $x_F\in\relint T_F$ for
      every nonempty proper $F$; and $X$ is a geometric realization of $W(\F(Q))$.
\item If $Q$ satisfies (I), then $Q$ is orthogonally reciprocal.
\item If $Q$ satisfies (I), then $X$ is an orthogonally bipolar polytope with poles
      $O,B$ (Definition~\ref{de:obp}), and all its vertices lie on the sphere with
      diameter $[O,B]$.
\end{enumerate}
\end{prop}
\begin{proo}
(a) We proceed in four steps.

\emph{The inequalities of $X$.} For a facet $\Phi$ of $Q$ let $\Phi^*$ be the dual vertex
of $Q^\circ$; thus $\Phi=\{y\in Q:\langle\Phi^*,y\rangle=1\}$, and $\Phi^*$ is the outer
normal of $\Phi$ scaled so that $\langle\Phi^*,y\rangle\le1$ on $Q$ (\cite[\S2.2]{Bipolar}).
Since a full-dimensional polytope is the intersection of the halfspaces of its facets,
\[
Q=\{y:\ \langle\Phi^*,y\rangle\le1\ \text{for every facet }\Phi\text{ of }Q\};
\]
equivalently, $Q=(Q^\circ)^\circ$, and it suffices to test the defining inequalities of
the polar at the vertices $\Phi^*$ of $Q^\circ$. A point
$x=(x_0,x')$ with $x_0>0$ lies in $T$ iff $x'/(\varepsilon x_0)\in Q$, i.e.\ iff
$\langle\Phi^*,x'\rangle\le\varepsilon x_0$ for all $\Phi$; for $x_0=0$ these inequalities force
$x'=O_d$, as the vectors $\Phi^*$ positively span $\R^d$. Hence $T$ is given by the
inequalities $\langle(-\varepsilon,\Phi^*),x\rangle\le0$, one for each facet $\Phi$, and the
equation of the corresponding facet $T_\Phi$ of $T$ is
$\langle(-\varepsilon,\Phi^*),x\rangle=0$. Since $\langle B,a_u\rangle=\beta$, the
inequality $\langle x-B,a_u\rangle\le0$ in the definition of $X$ reads
$\langle a_u,x\rangle\le\beta$. Altogether,
\[
X=\bigl\{x\in\R^{d+1}:\ \langle(-\varepsilon,\Phi^*),x\rangle\le0\ \text{for every facet }\Phi,\
\ \langle a_u,x\rangle\le\beta\ \text{for every vertex }u\bigr\},
\]
i.e.\ the points $x\in X$ are exactly those satisfying all these inequalities.

\emph{$X$ is bounded.} Every point of $T$ is $x=s(1,\varepsilon\eta)$ with $s\ge0$ and
$\eta\in Q$, and $\langle x-B,a_u\rangle=s\bigl(1+\varepsilon^2\langle\eta,u\rangle\bigr)-\beta$.
So $x\in X$ iff $s\bigl(1+\varepsilon^2h_Q(\eta)\bigr)\le\beta$, where
$h_Q(\eta)=\max_{u}\langle\eta,u\rangle$ is the support function. As $\eta\in Q$, we have
$h_Q(\eta)\ge\langle\eta,\eta\rangle\ge0$, hence $0\le s\le\beta$ on $X$. Thus $X$ is a
bounded polyhedron, i.e.\ a polytope.

\emph{The polar.} The point $z=(z_0,O_d)$ with $0<z_0<\beta$ satisfies all the
inequalities strictly, so it is an interior point of $X$. In the coordinates $y=x-z$
the inequalities become
\[
\langle(-\varepsilon,\Phi^*),y\rangle\le\varepsilon z_0\qquad\text{and}\qquad
\langle a_u,y\rangle\le\beta-z_0 ,
\]
with positive right-hand sides. If a polytope is given by finitely many inequalities,
$\{y\in\R^{d+1}:\langle n_i,y\rangle\le b_i\ \text{for } i=1,\dots,k\}$ with
$b_1,\dots,b_k>0$, then its polar is $\conv\{n_i/b_i: i=1,\dots,k\}$ \cite[Theorem~2.11(vii)]{Ziegler}. Therefore $(X-z)^\circ$ is the convex hull of the points
\[
\tfrac{1}{\varepsilon z_0}(-\varepsilon,\Phi^*)=\bigl(-\tfrac1{z_0},\ \tfrac1{\varepsilon z_0}\Phi^*\bigr)
\qquad\text{and}\qquad
\tfrac{1}{\beta-z_0}\,a_u=\bigl(\tfrac1{\beta-z_0},\ \tfrac{\varepsilon}{\beta-z_0}u\bigr).
\]
The first family consists of one point for each facet $\Phi$ of $Q$; all these points
have the $x_0$-coordinate $-1/z_0$, so they lie in the hyperplane $x_0=-1/z_0$. Their
$x'$-coordinates are $\frac1{\varepsilon z_0}\Phi^*$, and since $\Phi\mapsto\Phi^*$ is a
bijection from the facets of $Q$ onto the vertices of $Q^\circ$ (\cite[\S2.2]{Bipolar}),
these are exactly the vertices of the scaled polar $\frac1{\varepsilon z_0}Q^\circ$. Likewise
the second family lies in the hyperplane $x_0=1/(\beta-z_0)$, and its
$x'$-coordinates are exactly the vertices of $\frac{\varepsilon}{\beta-z_0}Q$. Hence
\[
(X-z)^\circ=\conv\Bigl(\bigl\{-\tfrac1{z_0}\bigr\}\times\tfrac1{\varepsilon z_0}Q^\circ\ \cup\
\bigl\{\tfrac1{\beta-z_0}\bigr\}\times\tfrac{\varepsilon}{\beta-z_0}Q\Bigr).
\]
The affine map $(x_0,x')\mapsto\bigl(\frac{x_0+1/z_0}{1/z_0+1/(\beta-z_0)},\,x'\bigr)$
sends the two hyperplanes to the heights $0$ and $1$; so $(X-z)^\circ$ is affinely
equivalent to the Cayley polytope
$\mathrm{Cay}\bigl(\frac1{\varepsilon z_0}Q^\circ,\frac{\varepsilon}{\beta-z_0}Q\bigr)$. (The order of
the two layers is immaterial: $x_0\mapsto1-x_0$ exchanges them.)

\emph{The face lattice.} By \cite[\S2.2]{Bipolar}, the face lattice of a Cayley polytope
does not change under positive scaling of its two layers, so
$\F((X-z)^\circ)\cong\F(\mathrm{Cay}(Q,Q^\circ))$, and polarity gives
$\F(X)\cong\F(\mathrm{Cay}(Q,Q^\circ))^{\mathrm{op}}$. Finally, every vertex of either
layer is a vertex of the Cayley polytope, so each of the points $n_i/b_i$ above is a
vertex of $(X-z)^\circ$. Hence none of the inequalities describing $X$, namely
\[
\langle(-\varepsilon,\Phi^*),x\rangle\le0\quad(\Phi\text{ a facet of }Q)\qquad\text{and}\qquad
\langle a_u,x\rangle\le\beta\quad(u\text{ a vertex of }Q),
\]
is redundant, and each of them defines a facet of $X$: the inequality for $\Phi$
defines the facet $X\cap\aff T_\Phi$ (the equation of $\aff T_\Phi$ was found in the
first step, ``The inequalities of $X$''), and the inequality for $u$ the facet
$X\cap H_u$.

(b) By \cite[Remark~5.3]{Bipolar}, the vertex of $X$ dual to the facet
of the Cayley polytope labelled $[F,F]$ lies exactly on the facets $H_u$, $u\in F$, and
$\aff T_\Phi$, $\Phi\supseteq F$. The hyperplanes $\aff T_\Phi$, $\Phi\supseteq F$,
meet in the linear span of $T_F$, and on that span the conditions
$\langle x-B,a_u\rangle=0$, $u\in F$, single out $x_F$. For $F=\emptyset$ and $F=Q$
one obtains $O$ and $B$. A point of $T$ that lies on the facets of $T$ containing $T_F$
and on no other facet of $T$ is in $\relint T_F$. The facet structure is that of
$W(\F(Q))$ by \cite[Theorem~5.2]{Bipolar}.

(c) Let $F$ be a nonempty proper face of $Q$; we verify ($\star1$) and ($\star2$)
of Definition~\ref{def:star}.

($\star1$): By (b), $x_F\in\relint T_F$. By Lemma~\ref{lem:tall}(b) this is equivalent
to $q_F\in\relint F$.

($\star2$): Let $u$ be a vertex of $Q$ with $u\notin F$. By (b), the vertices of $X$ on
its facet $X\cap H_u$ are the points $x_G$ with $u\in G$ (the facet ${\uparrow}u$ of
$W(\F(Q))$); hence $x_F\notin H_u$. Since $x_F\in X$, this gives
$\langle x_F-B,a_u\rangle<0$, and by the formula of Lemma~\ref{lem:tall}, whose factor
$\varepsilon^2\beta/(1+\varepsilon^2m)$ is positive, $\langle q_F,u\rangle<m=\|q_F\|^2$.

(d) By (b), $X$ is a geometric bipolar polytope with poles $x_\emptyset=O$ and
$x_Q=B$; the atoms of its vertex lattice are the points $x_{\{u\}}$ on the edges of
$T$. The face $[O,x_{\{u\}}]$ is a segment in the direction $a_u$, and the face
$[x_{\{u\}},B]$ is the facet $X\cap H_u$, orthogonal to $a_u$; so they are
perpendicular, which is the defining condition of an orthogonally bipolar polytope.
Every vertex $x_F$ satisfies $\langle x_F-O,x_F-B\rangle=0$, because $x_F-B$ is
orthogonal to the span of $T_F$, which contains $x_F$; this is the Thales condition
for the sphere with diameter $[O,B]$ (see also Theorem~\ref{theo:obpthales}).
\end{proo}

Conversely, every orthogonally bipolar realization arises in this way, so (I) is also
necessary up to the choice of the realization.

\begin{theo}[orthogonally bipolar realizations]\label{prop:obpreal}
For a polytope $Q$ the following are equivalent:
\begin{enumerate}[label=(\roman*)]
\item $W(\F(Q))$ has a geometric realization that is an orthogonally bipolar polytope;
\item some realization $Q'$ of the combinatorial type of $Q$, with $O\in\inte Q'$,
      satisfies (I), i.e.\ is perfectly centered (Theorem~\ref{conj:I}).
\end{enumerate}
Moreover, every orthogonally bipolar realization of $W(\F(Q))$ is, up to an isometry,
the apex cutout $X$ of Proposition~\ref{prop:cutout} for such a $Q'$, with $\varepsilon=1$.
In particular, every orthogonally bipolar realization of $W(\F(Q))$ is inscribed in the
sphere whose diameter joins its poles (Proposition~\ref{prop:cutout}(d)).
\end{theo}
\begin{proo}
(ii)$\Rightarrow$(i) is Proposition~\ref{prop:cutout}(b),(d) applied to $Q'$.

(i)$\Rightarrow$(ii). Let $P$ be an orthogonally bipolar polytope realizing $W(\F(Q))$,
with poles $\sigma,\nu$. After an isometry, $\sigma=O$ and $\nu=(\beta,O_d)$ with
$\beta>0$. Let $T$ be the tangent cone of $P$ at $O$.

\emph{The base $Q'$.} Every atom $a$ lies on the sphere with diameter $[O,\nu]$
Theorem~\ref{theo:obpthales}, so $\langle a,\nu\rangle=\|a\|^2>0$; hence all edges of $T$
point into the halfspace $x_0>0$, and $T\cap\{x_0=1\}=\{1\}\times Q'$ for a polytope
$Q'$. The pole $\nu$ lies in no facet of $P$ containing $O$, hence in the interior of
$T$, so $O\in\inte Q'$. The faces of $T$ are the cones $T_F$ over the faces $F$ of $Q'$,
and $\F(Q')\cong\F(P/\sigma)\cong\F(Q)$ by \cite[Proposition~3.13]{Bipolar}; so $Q'$ is a
realization of the type of $Q$. From now on we use the notation introduced before
Lemma~\ref{lem:tall} for $Q'$ in place of $Q$, with $\varepsilon=1$ and $B=\nu$: thus
$a_u=(1,u)$ for the vertices $u$ of $Q'$, so that $T=\cone\{a_u\}$ is the tangent cone
above; $x_F$ is the orthogonal projection of $B$ onto $\operatorname{lin}T_F$;
$H_u=\{x:\langle x-B,a_u\rangle=0\}$; and
$X=T\cap\bigcap_u\{x:\langle x-B,a_u\rangle\le0\}$.

\emph{The vertices.} Let $x$ be a vertex of $P$ and $F$ the face of $Q'$ with
$\aff[\sigma,x]=\operatorname{lin}T_F$. By complementarity
Lemma~\ref{lm:full}, $[\sigma,x]\perp[x,\nu]$, so $\nu-x$ is orthogonal to
$\operatorname{lin}T_F$; as $x\in\operatorname{lin}T_F$, the point $x$ is the orthogonal
projection of $B=\nu$ onto $\operatorname{lin}T_F$, i.e.\ $x=x_F$. So the vertices of
$P$ are exactly the points $x_F$.

\emph{The facets of $P$.} We show that $P=X$. The facets of $P$ containing $\sigma$
lie in the hyperplanes $\aff T_\Phi$, $\Phi$ a facet of $Q'$, which bound $T$. A facet
$[a,\nu]$ with $a=x_{\{u\}}$ an atom has the vertices $x_G$, $u\in G$, which lie in
$H_u$ (Lemma~\ref{lem:tall}(a)); so it lies in $H_u$. A polytope is the
intersection of the closed halfspaces that are bounded by the hyperplanes of its
facets and contain it. For the facets in $\aff T_\Phi$ these halfspaces are those
bounding $T$. For a facet in $H_u$ it is one of the two halfspaces
$\{x:\langle x-B,a_u\rangle\le0\}$ and $\{x:\langle x-B,a_u\rangle\ge0\}$; since $O\in P$
and $\langle O-B,a_u\rangle=-\langle B,a_u\rangle=-\beta<0$, it is the first one. Hence
$P=T\cap\bigcap_u\{x:\langle x-B,a_u\rangle\le0\}=X$. Moreover
$x_F\in H_u$ iff $u\in F$, and $x_F\in\aff T_\Phi$ iff $F\subseteq\Phi$, since these are
the facet incidences of $W(\F(Q'))$.

\emph{(I) for $Q'$.} Let $C=\mathrm{Cay}(Q',Q'^\circ)$. We call a face of $C$
\emph{mixed} if it has vertices in both layers. By the Cayley trick
(\cite[\S2.2]{Bipolar}), the mixed faces of $C$ correspond to the pairs of faces of $Q'$
and $Q'^\circ$ with a common normal in their relative interiors; writing the face of
$Q'^\circ$ as $G^*$ and using $N_{Q'^\circ}(G^*)=\cone(G)$, these are the pairs
$(H,G)$ with $\relint N_{Q'}(H)\cap\relint\cone(G)\ne\emptyset$, the mixed face having
the vertices of $H$ in the lower layer and those of $G^*$ in the upper one. So (I) says
that the mixed faces of $C$ are exactly the faces with lower part $H$ and upper part
$G^*$ for the pairs $H\subseteq G$. We verify this in two steps.

First, the facets of $C$. By Proposition~\ref{prop:cutout}(a), $\F(X)$ is
anti-isomorphic to $\F(C)$: vertices of $X$ correspond to facets of $C$, facets of $X$
to vertices of $C$, incidences are preserved, the facet $X\cap H_u$ corresponds to the
lower vertex given by $u$, and the facet on $\aff T_\Phi$ to the upper vertex given by
$\Phi^*$. By the previous step, $x_F$ lies on $H_u$ iff $u\in F$, and on
$\aff T_\Phi$ iff $F\subseteq\Phi$. Hence the facet $C_F$ of $C$ dual to the vertex
$x_F$, for a nonempty proper face $F$, has the vertices of $F$ in the lower layer and
those of $F^*$ in the upper one.

Second, the mixed faces. Every face of $C$ is the intersection of the facets
containing it, and a mixed face lies in no base, so it is $\bigcap_iC_{F_i}$ for some
nonempty proper faces $F_1,\dots,F_k$. Its lower part is $H=\bigcap_iF_i$ and its upper
part is $\bigcap_iF_i^*=G^*$, where $G$ is the smallest face containing all $F_i$; so
$H\subseteq G$. Conversely, for $\emptyset\ne H\subseteq G\ne Q'$ the facets $C_K$ with
$H\subseteq K\subseteq G$ intersect in a face with lower part $H$ and upper part $G^*$,
both nonempty; so it is a mixed face. Thus the mixed faces are exactly those given by
the pairs $H\subseteq G$, which is (I).
\end{proo}

Under (I), the abstract construction $W(\F(Q))$, the Cayley realization and the apex
cutout by orthogonal projections thus all coincide combinatorially, and by
Theorem~\ref{prop:obpreal}, $W(\F(Q))$ has an orthogonally bipolar realization exactly when the
combinatorial type of $Q$ has a realization satisfying (I). Proposition~\ref{prop:cutout}(c)
shows that (I) implies orthogonal reciprocity; the converse holds as well, by Broadie's
theorem (Theorem~\ref{conj:I} below).

The equivalence of (i) and (iii) below is, in the language of balanced pairs, Dobbins'
Lemma~3.3 \cite{Dobbins}, derived from Broadie's theorem \cite{Broadie}; condition (ii) is
the definition of perfect centering by Fukuda and Weibel \cite{FukudaWeibel}. These results, together
with the Cayley trick, are collected in \cite[Theorem~5.2]{Bipolar}, which we use in the
proof.

\begin{theo}[perfect centering]\label{conj:I}
For a polytope $Q$ with $O\in\inte Q$ the following are equivalent:
\begin{enumerate}[label=(\roman*)]
\item $Q$ satisfies the interval condition (I);
\item $Q$ is orthogonally reciprocal;
\item condition ($\star1$) holds for every nonempty proper face of $Q$.
\end{enumerate}
\end{theo}
\begin{proo}
(i)$\Rightarrow$(ii) is Proposition~\ref{prop:cutout}(c), (ii)$\Rightarrow$(iii) is trivial, and
(iii)$\Rightarrow$(i) is \cite[Theorem~5.2]{Bipolar}, (e)$\Rightarrow$(a).
\end{proo}

In particular, condition ($\star1$) for all faces implies ($\star2$). This holds only for all
faces together: for a vertex $v$, condition ($\star1$) is automatic since $q_v=v$, whereas
($\star2$) requires that $v$ is the unique maximizer of $\langle v,\cdot\rangle$ on $Q$, which can
fail for a single vertex. In Broadie's terms, the projection condition (PR) for all faces
is equivalent to the supporting hyperplane condition (SH) for all faces
\cite[Lemma~4.4]{Broadie}.

Condition (I) depends on the position of the origin $O$ inside $Q$. By
Theorem~\ref{conj:I}, it holds about $O$ exactly when the foot of the perpendicular from $O$
to the affine hull of every face lies in the relative interior of that face. This leads to
the following notion, which we use for polytopes of any dimension.

\begin{de}[orthogonal kernel]\label{def:pedal}
Let $P\subset\R^n$ be an $n$-polytope. For $B\in\R^n$ and a nonempty face $G$ of $P$ let
$x_{B,G}$ be the orthogonal projection of $B$ onto $\aff G$. The \emph{orthogonal kernel} of
$P$ is
\[
\ker_\perp(P)=\bigl\{B\in\inte P:\ x_{B,G}\in\relint G\ \text{ for every nonempty face } G
\text{ of } P\bigr\}.
\]
\end{de}

For polygons, $\ker_\perp(P)$ is Gr\"unbaum's polarity kernel
\cite[Proposition~2]{GrunbaumAnti2}. By definition,
$B\in\ker_\perp(P)$ if and only if $P$ is perfectly centered about $B$ in the sense of
Broadie \cite{Broadie}, Fukuda and Weibel \cite{FukudaWeibel} and Dobbins \cite{Dobbins};
Lawrence \cite{Lawrence} calls such polytopes \emph{adequately centered}.

\begin{prop}[description of the orthogonal kernel]\label{rem:centres}
For a polytope $Q$, the orthogonal kernel $\ker_\perp(Q)$ is the set of points $O\in\inte Q$
about which $Q$ satisfies (I). Moreover
\[
\ker_\perp(Q)=\bigcap_{F}\bigl(\relint F-\relint N_Q(F)\bigr),
\]
the intersection, over the nonempty proper faces $F$ of $Q$, of the open ``inner normal
prisms'' of the faces. In particular, $\ker_\perp(Q)$ is the interior of a convex polytope, possibly
empty.
\end{prop}
\begin{proo}
For $O\in\inte Q$ and a nonempty proper face $F$, let $q_F$ be the foot of $O$ on $\aff F$
and $p=q_F-O$. Condition ($\star1$) says $q_F\in\relint F$, and ($\star2$) says that
$\langle p,\cdot\rangle$ attains its maximum on $Q$ exactly on $F$, i.e.\
$p\in\relint N_Q(F)$. We claim that $O$ satisfies both conditions for $F$ if and only if
$O=f-n$ for some $f\in\relint F$ and some $n\in\relint N_Q(F)$. If both conditions hold, take
$f=q_F$ and $n=p$; then $f-n=q_F-(q_F-O)=O$. Conversely, let $O=f-n$ with $f$ and $n$ as
stated. Every $n\in N_Q(F)$ is orthogonal to the direction space of $\aff F$, since
$\langle n,\cdot\rangle$ attains its maximum on all of $F$ and is therefore constant on $F$.
So $O-f=-n$ is orthogonal to $\aff F$ and $f\in\aff F$, which means that $f$ is the foot of
the perpendicular from $O$ to $\aff F$, i.e.\ $f=q_F$. Hence $q_F\in\relint F$, which is
($\star1$), and $p=q_F-O=n\in\relint N_Q(F)$, which is ($\star2$). By Theorem~\ref{conj:I}, $O\in\ker_\perp(Q)$ if and only if
this holds for every $F$. Each set $\relint F-\relint N_Q(F)$ is an open convex
polyhedral set, and for a facet $F$ it lies on the inner side of $\aff F$; so the
intersection is a bounded open convex polyhedral set, the interior of a convex polytope
if it is nonempty.
\end{proo}

\begin{figure}[htbp]
\centering
\begin{tikzpicture}[scale=0.62,line cap=round,line join=round]
\fill[gray!8] (0.000,0.000)--(5.000,0.000)--(6.000,2.200)--(2.200,4.000)--(-0.800,1.600)--cycle;
\draw[gray!55,densely dashed,line width=0.35pt] (0.000,0.000)--(0.000,2.240);
\draw[gray!55,densely dashed,line width=0.35pt] (5.000,0.000)--(5.000,2.674);
\draw[gray!55,densely dashed,line width=0.35pt] (5.000,0.000)--(0.026,2.261);
\draw[gray!55,densely dashed,line width=0.35pt] (6.000,2.200)--(2.142,3.954);
\draw[gray!55,densely dashed,line width=0.35pt] (6.000,2.200)--(4.958,0.000);
\draw[gray!55,densely dashed,line width=0.35pt] (2.200,4.000)--(0.305,0.000);
\draw[gray!55,densely dashed,line width=0.35pt] (2.200,4.000)--(5.145,0.319);
\draw[gray!55,densely dashed,line width=0.35pt] (-0.800,1.600)--(0.480,0.000);
\draw[gray!55,densely dashed,line width=0.35pt] (-0.800,1.600)--(3.124,3.562);
\draw[gray!55,densely dashed,line width=0.35pt] (0.000,0.000)--(5.178,2.589);
\fill[blue!25] (3.857,1.929)--(2.714,3.357)--(1.641,2.821)--(1.137,1.756)--(2.381,1.190)--cycle;
\draw[blue!70!black,line width=0.7pt] (3.857,1.929)--(2.714,3.357)--(1.641,2.821)--(1.137,1.756)--(2.381,1.190)--cycle;
\draw[line width=0.9pt] (0.000,0.000)--(5.000,0.000)--(6.000,2.200)--(2.200,4.000)--(-0.800,1.600)--cycle;
\node[blue!60!black,font=\small] at (2.346,2.211) {$\ker_\perp$};
\node[font=\small] at (2.6,-0.75) {(a)};
\end{tikzpicture}\hspace{0.4cm}
\begin{tikzpicture}[scale=0.62,line cap=round,line join=round]
\fill[gray!8] (0.000,0.000)--(5.500,0.000)--(4.200,2.200)--(0.600,3.600)--(-1.400,1.000)--cycle;
\fill[orange!35] (0.477,-0.000)--(5.500,-0.000)--(4.200,2.200)--cycle;
\draw[gray!70,densely dashed,line width=0.4pt] (4.200,2.200)--(0.477,0.000);
\draw[black,line width=0.45pt] (4.363,1.925)--(4.087,1.762)--(3.925,2.037);
\fill[green!30] (-1.400,1.000)--(0.000,0.000)--(2.143,3.000)--(0.600,3.600)--cycle;
\draw[gray!70,densely dashed,line width=0.4pt] (0.000,0.000)--(2.143,3.000);
\draw[black,line width=0.45pt] (-0.260,0.186)--(-0.074,0.446)--(0.186,0.260);
\draw[line width=0.9pt] (0.000,0.000)--(5.500,0.000)--(4.200,2.200)--(0.600,3.600)--(-1.400,1.000)--cycle;
\draw[orange!80!black,line width=1.6pt] (5.500,0.000)--(4.200,2.200);
\draw[green!50!black,line width=1.6pt] (-1.400,1.000)--(0.000,0.000);
\node[font=\small] at (2.1,-0.75) {(b)};
\end{tikzpicture}\hspace{0.4cm}
\begin{tikzpicture}[scale=0.9,line cap=round,line join=round]
\filldraw[fill=blue!22,draw=blue!60!black,line width=0.35pt] (0.594,-0.788)--(0.167,-0.755)--(-0.371,-0.870)--(0.000,-0.976)--cycle;
\filldraw[fill=blue!22,draw=blue!60!black,line width=0.35pt] (-0.723,-0.550)--(-0.371,-0.870)--(-0.556,-0.329)--(-0.965,-0.082)--cycle;
\filldraw[fill=blue!22,draw=blue!60!black,line width=0.35pt] (-0.723,0.426)--(-0.594,0.788)--(-0.890,0.206)--(-0.965,-0.082)--cycle;
\filldraw[fill=blue!22,draw=blue!60!black,line width=0.35pt] (-0.556,-0.329)--(-0.371,-0.870)--(0.167,-0.755)--(0.223,-0.356)--cycle;
\filldraw[fill=blue!22,draw=blue!60!black,line width=0.35pt] (0.223,-0.356)--(0.167,-0.755)--(0.594,-0.788)--(0.890,-0.206)--cycle;
\filldraw[fill=blue!22,draw=blue!60!black,line width=0.35pt] (0.890,-0.206)--(0.594,0.513)--(0.167,0.221)--(0.223,-0.356)--cycle;
\filldraw[fill=blue!22,draw=blue!60!black,line width=0.35pt] (-0.556,-0.329)--(-0.371,0.431)--(-0.723,0.426)--(-0.965,-0.082)--cycle;
\filldraw[fill=blue!22,draw=blue!60!black,line width=0.35pt] (0.223,-0.356)--(0.167,0.221)--(-0.371,0.431)--(-0.556,-0.329)--cycle;
\filldraw[fill=blue!22,draw=blue!60!black,line width=0.35pt] (-0.594,0.788)--(-0.723,0.426)--(-0.371,0.431)--(0.000,0.976)--cycle;
\filldraw[fill=blue!22,draw=blue!60!black,line width=0.35pt] (0.000,0.976)--(-0.371,0.431)--(0.167,0.221)--(0.594,0.513)--cycle;
\filldraw[fill=blue!22,draw=blue!60!black,line width=0.35pt] (0.965,0.082)--(0.723,0.550)--(0.594,0.513)--(0.890,-0.206)--cycle;
\filldraw[fill=blue!22,draw=blue!60!black,line width=0.35pt] (0.371,0.870)--(0.000,0.976)--(0.594,0.513)--(0.723,0.550)--cycle;
\draw[line width=0.9pt] (1.113,0.658)--(0.000,1.951);
\draw[gray!70,dashed,line width=0.45pt] (-1.781,0.411)--(1.113,0.658);
\draw[line width=0.9pt] (1.781,-0.411)--(0.000,1.951);
\draw[line width=0.9pt] (-1.113,-0.658)--(0.000,1.951);
\draw[gray!70,dashed,line width=0.45pt] (-1.781,0.411)--(0.000,-1.951);
\draw[line width=0.9pt] (1.781,-0.411)--(-1.113,-0.658);
\draw[line width=0.9pt] (-1.781,0.411)--(0.000,1.951);
\draw[line width=0.9pt] (1.781,-0.411)--(1.113,0.658);
\draw[line width=0.9pt] (1.781,-0.411)--(0.000,-1.951);
\draw[gray!70,dashed,line width=0.45pt] (1.113,0.658)--(0.000,-1.951);
\draw[line width=0.9pt] (-1.781,0.411)--(-1.113,-0.658);
\draw[line width=0.9pt] (-1.113,-0.658)--(0.000,-1.951);
\node[blue!60!black,font=\small] at (0.000,0.000) {$\ker_\perp$};
\node[font=\small] at (0,-2.45) {(c)};
\end{tikzpicture}
\caption{The orthogonal kernel (Proposition~\ref{rem:centres}). (a) For a polygon, the orthogonal
kernel $\ker_\perp$ is the intersection of the strips perpendicular to the edges (their boundaries are
dashed); here it is a pentagon covering about $18\%$ of the area. (b) For this
pentagon the strips of the two marked edges (right angles marked) are already
disjoint, so $\ker_\perp=\emptyset$:
no choice of the origin makes this realization satisfy (I). (c) For the regular
octahedron $\conv\{\pm e_i\}$, the orthogonal kernel is the polytope with the $26$
vertices $(\pm\frac12,0,0)$, $(\pm\frac13,\pm\frac13,0)$, $(\pm\frac14,\pm\frac14,\pm\frac14)$
and their permutations; it is cut out by the planes perpendicular to the edges at
their endpoints and the planes through the edges perpendicular to the facets.}
\label{fig:centres}
\end{figure}
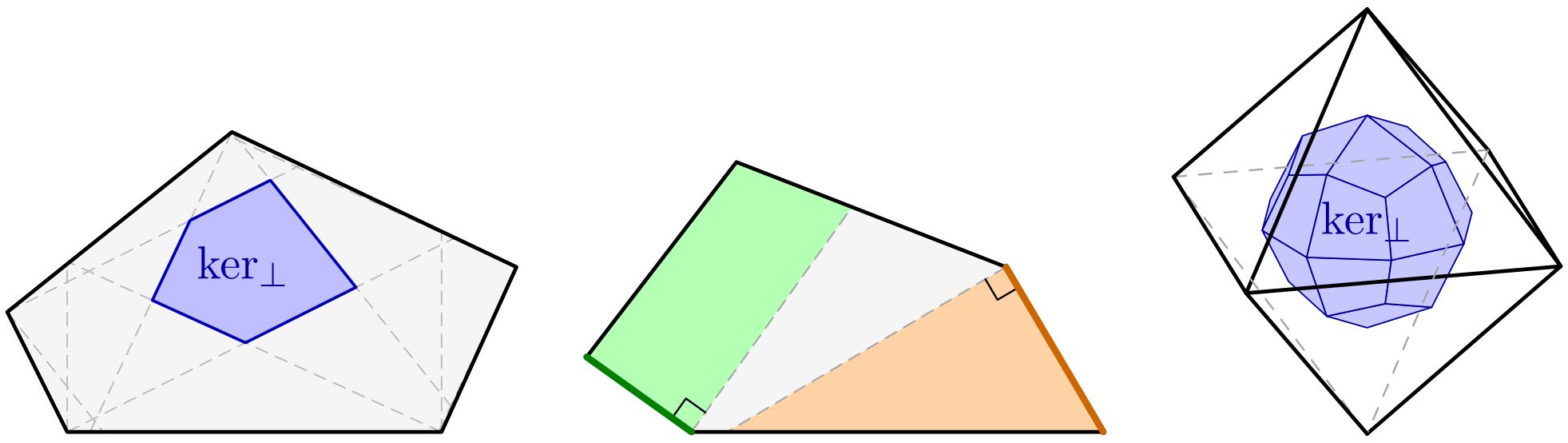

For a polygon, $\ker_\perp(Q)$ is the intersection of the open strips perpendicular to the edges;
this is Gr\"unbaum's description of the polarity kernel of a polygon
\cite[Proposition~2]{GrunbaumAnti2}, which may be empty, as for a suitable parallelogram
\cite[Figure~5]{GrunbaumAnti2}. Proposition~\ref{rem:centres} extends it to polytopes of
every dimension.
Numerically, $\ker_\perp(Q)$ is the whole interior for the cube and the regular tetrahedron and
about a fifth of it for the octahedron. For a fixed realization, $\ker_\perp(Q)\neq\emptyset$ is
the feasibility of finitely many strict linear inequalities in $O$. Some combinatorial
$4$-polytopes have no realization with $\ker_\perp(Q)\neq\emptyset$: Dobbins
\cite[Theorem~3.9]{Dobbins} constructed one without any antiprism, whereas by
Theorem~\ref{conj:I} a realization with $\ker_\perp(Q)\neq\emptyset$ would make $\mathrm{Cay}(Q,Q^\circ)$
an antiprism. Figure~\ref{fig:centres} shows the region for two pentagons, one with
$\ker_\perp=\emptyset$, and for the regular octahedron.

Perfect centering is inherited by the interval polytope, by a theorem of Lawrence
\cite[Theorem~2.2]{Lawrence}. Together with Theorem~\ref{prop:obpreal}, this gives
orthogonally bipolar realizations of all iterated interval polytopes.

\begin{coro}[iterated interval polytopes]\label{cor:iterate}
If the combinatorial type of $Q$ has a perfectly centered realization, then so does the
interval polytope $W(\F(Q))$. Consequently, for every $k\ge1$, the $k$-th iterated interval
polytope of $Q$, whose face lattice is obtained from $\F(Q)$ by taking the lattice of
intervals $k$ times, has a realization, and it has an orthogonally bipolar realization
inscribed in a sphere.
\end{coro}
\begin{proo}
Let $Q$ be perfectly centered about $O$. By Lawrence \cite[Theorem~2.2]{Lawrence}, the
antiprism $A=\conv\bigl((Q\times\{-1\})\cup(Q^\circ\times\{1\})\bigr)$ is perfectly centered
about the origin, and so is its polar $A^\circ$, since a polytope is perfectly centered if
and only if its polar is \cite[Lemma~4.4]{Broadie}. By Broadie's theorem, $A^\circ$
realizes $W(\F(Q))$ \cite[Theorem~5.2]{Bipolar}, \cite[Lemma~2.1]{Lawrence}. This proves
the first statement, and the second follows by induction, together with
Theorem~\ref{prop:obpreal}.
\end{proo}

Lawrence's binary partition polytopes \cite{Lawrence} arise in this way from a segment; by
Corollary~\ref{cor:iterate}, each of them has an orthogonally bipolar realization inscribed
in a sphere.

Both the following application to bases and the tilings of Section~\ref{sec:pedal} rest
on one statement about cutouts at the apex of a cone.

\begin{prop}[orthogonal cutouts at the apex of a cone]\label{prop:conecut}
Let $T\subset\R^n$ be a pointed $n$-dimensional polyhedral cone with apex $v$, let
$B\in\inte T$, and for each face $F$ of $T$ let $x_F$ be the orthogonal projection of $B$ onto
$\aff F$. If $x_F\in\relint F$ for every face $F\neq\{v\}$ of $T$, then $\conv\{x_F\}$ is an
orthogonally bipolar polytope with poles $v,B$, whose vertices are the points $x_F$ and whose
vertex lattice is the face lattice of $T$.
\end{prop}
\begin{proo}
Take $v$ as the origin, so that $\aff F$ is the linear span $L_F$ of $F$. For an edge
$e=\R_{\ge0}u$ of $T$ the hypothesis $x_e\in\relint e$ gives $\langle u,B\rangle>0$; as every
nonzero vector of $T$ is a nonnegative combination of edge directions, $\langle y,B\rangle>0$
for all $y\in T\setminus\{0\}$. Hence the hyperplane $H=\{y:\langle y,B\rangle=\|B\|^2\}$ through
$B$ perpendicular to $vB$ meets every ray of $T$ in one point, and $Q'=T\cap H$ is a polytope
whose faces are the sets $F'=F\cap H$, $F\neq\{v\}$. For $y\in L_F$ we have
$\langle y,B\rangle=\langle y,x_F\rangle$, since $B-x_F$ is orthogonal to $L_F$; so
$\aff F'=\{y\in L_F:\ \langle y,x_F\rangle=\|B\|^2\}$, and for $y\in H$ we have
$\|y-B\|^2=\|y\|^2-\|B\|^2$ by Pythagoras, since $y-B\perp B$. Hence the foot of the
perpendicular from $B$ to $\aff F'$ is the point of $\aff F'$ nearest to $v$, which is the
positive multiple $(\|B\|^2/\|x_F\|^2)\,x_F$ of $x_F$; it lies in $\relint F'$ if and only if
$x_F\in\relint F$. So $Q'$ is perfectly centered about $B$, by Theorem~\ref{conj:I} it
satisfies the interval condition about $B$, and Corollary~\ref{cor:lines}(b), applied to the
lines carrying the edges of $T$ and to the poles $v,B$, gives the claim.
\end{proo}

The simplest cutouts are those at the apex of a simplicial cone. They connect the theory
with orthogonal and non-orthogonal expansions. Let $b_1,\dots,b_n$ be a basis of $\R^n$ and
$\nu\in\R^n$. The basis spans the simplicial cone $T=\cone\{b_1,\dots,b_n\}$ with apex $O$;
cutting it off by a hyperplane that meets all its rays gives a pyramid $P$ with apex $O$,
whose faces containing $O$ are the cones $T_S=\cone\{b_i:i\in S\}$, $S\subseteq\{1,\dots,n\}$,
cut off, with the affine hulls $\operatorname{span}\{b_i:i\in S\}$. For each $S$ let $x_S$ be
the orthogonal projection of $\nu$ onto $\operatorname{span}\{b_i:i\in S\}$, i.e.\ the best
approximation of $\nu$ by the vectors $b_i$, $i\in S$; so $x_\emptyset=O$ and
$x_{\{1,\dots,n\}}=\nu$, and we write $x_S=\sum_{i\in S}y_{S,i}\,b_i$. By
Lemma~\ref{lm:orth-thales}, applied to $P$, its vertex $O$ and the point $\nu$, all points $x_S$
lie on the sphere with diameter $[O,\nu]$; note that the lemma does not require
$\nu\in P$, and that the points $x_S$ do not depend on the choice of the truncation. The
question is whether the $x_S$ are the vertices of an orthogonal cutout, i.e.\ of an
orthogonally bipolar polytope with poles $O$ and $\nu$. Its vertex lattice is then the
Boolean lattice of the faces of $P$ containing $O$, so it is a Hjelmslev cuboid, and the
sphere with diameter $[O,\nu]$ is its circumscribed sphere (Theorem~\ref{theo:obpthales}).

\begin{prop}[Hjelmslev cuboids from bases]\label{prop:bases}
The polytope $\conv\{x_S\}$ is a Hjelmslev $n$-cuboid with poles $O,\nu$ whose vertex
lattice is the Boolean lattice of the sets $S$ if and only if there are signs
$\varepsilon_1,\dots,\varepsilon_n\in\{1,-1\}$ such that $\varepsilon_i\,y_{S,i}>0$ for all
$S$ and all $i\in S$; that is, the coefficient of $b_i$ in the best approximation of $\nu$
by any subfamily containing $b_i$ never changes its sign.
\end{prop}
\begin{proo}
The points $x_S$ do not change if $b_i$ is replaced by $-b_i$, so we may assume
$\varepsilon_i=1$ in the condition. Suppose it holds. Then the condition says that $x_S\in\relint T_S$ for every
nonempty $S$, where $x_S$ is the orthogonal projection of $\nu$ onto $\aff T_S$; in
particular $\nu=x_{\{1,\dots,n\}}\in\inte T$. By Proposition~\ref{prop:conecut}, $\conv\{x_S\}$
is an orthogonally bipolar polytope with poles $O,\nu$ whose vertex lattice is the face
lattice of the simplicial cone $T$, i.e.\ the Boolean lattice of the sets $S$. Conversely, let $\conv\{x_S\}$
be such a Hjelmslev cuboid. Its atoms are the points $x_{\{i\}}=y_{\{i\},i}b_i$, so its
tangent cone at $O$ is $\cone\{\varepsilon_ib_i\}$ with $\varepsilon_i$ the sign of
$y_{\{i\},i}$, and by Proposition~\ref{prop:obp-cut} every vertex $x_S$ lies in the relative
interior of the face $\cone\{\varepsilon_ib_i:i\in S\}$ of that cone, i.e.\
$\varepsilon_iy_{S,i}>0$ for $i\in S$.
\end{proo}

\begin{rem}[orthonormal bases and Fourier series]\label{rem:bases}
For an orthonormal basis, $y_{S,i}=\langle\nu,b_i\rangle$ does not depend on $S$, so the
condition of Proposition~\ref{prop:bases} holds whenever all $\langle\nu,b_i\rangle\neq0$, and
$\conv\{x_S\}$ is the box with opposite vertices $O$ and $\nu$; its vertices are the partial
sums $x_S=\sum_{i\in S}\langle\nu,b_i\rangle b_i$ of the expansion of $\nu$. The same holds in a
Hilbert space. For an orthonormal system $(e_k)$ let $S_k$ denote the orthogonal projection
onto $\operatorname{span}\{e_1,\dots,e_k\}$. For a vector $f$, the partial sums
$S_kf=\sum_{i\le k}\langle f,e_i\rangle e_i$ of the Fourier series of $f$ are thus orthogonal
projections of $f$, so they satisfy the Thales condition $\langle S_kf,f-S_kf\rangle=0$
(cf.\ Lemma~\ref{lm:orth-thales}) and lie on the sphere with diameter $[0,f]$. Written out,
the Thales condition is the identity
$\|S_kf\|^2=\langle S_kf,f\rangle=\sum_{i\le k}|\langle f,e_i\rangle|^2$; Bessel's inequality
$\sum_i|\langle f,e_i\rangle|^2\le\|f\|^2$ expresses that every point of this sphere has norm at
most $\|f\|$, and Parseval's identity expresses that for a complete system the partial sums
converge to the antipode $f$ of $0$. The sequence $0=S_0f,S_1f,S_2f,\dots$ runs along edges of
these boxes. Passing to the limit, let $c_i=\langle f,e_i\rangle$ and, replacing $e_i$ by
$-e_i$ if necessary, assume $c_i\ge0$. Since $c_i\to0$, the set
$\{\sum_it_ic_ie_i:\ 0\le t_i\le1\}$ is a Hilbert cube, a compact convex set in which the
convex hulls of these boxes are dense; its extreme points are the corners
$x_A=\sum_{i\in A}c_ie_i$ for all sets $A$ of indices, finite or infinite. Each $x_A$ is the
orthogonal projection of $f$ onto the closed span of $\{e_i:i\in A\}$, so all corners lie on
the sphere with diameter $[0,f]$: this sphere is circumscribed about the Hilbert cube, and the
Fourier series of $f$ is an edge path from the corner $0$ to the opposite corner $f$.

For a non-orthogonal basis the coefficients $y_S=G_S^{-1}(\langle b_i,\nu\rangle)_{i\in S}$
depend on $S$, where $G_S=(\langle b_i,b_j\rangle)_{i,j\in S}$ is the Gram matrix of the
subfamily, and Proposition~\ref{prop:bases} says that the picture remains convex, a Hjelmslev
cuboid instead of a box, exactly when no coefficient changes its sign. This holds, for
instance, if the basis is obtuse, $\langle b_i,b_j\rangle\le0$ for $i\neq j$, and
$\langle b_i,\nu\rangle>0$ for all $i$: the Gram matrix is then a Stieltjes matrix, and the
inverses of its principal submatrices are entrywise nonnegative with positive diagonal
\cite{BermanPlemmons}, so that all $y_{S,i}>0$. In general the condition may fail.
\end{rem}

\section{Covering a polytope by Thales balls}\label{sec:app}

Section~\ref{sec:cutouts} concerned a single vertex $V$: the orthogonal projections of a
point $B$ onto the faces through $V$ lie on the sphere with diameter $[V,B]$
(Lemma~\ref{lm:orth-thales}), and when they form an admissible cutout, this is an
orthogonally bipolar polytope inscribed in that sphere (Proposition~\ref{prop:cut-obp} and
Theorem~\ref{theo:obpthales}). We now consider all vertices of a polytope together: the balls
with diameters $[V,B]$ cover the polytope (Theorem~\ref{thm:cover}), and the power diagram of
their boundary spheres refines this covering to a tiling (Corollary~\ref{cor:power}).

Throughout this section $P\subseteq\R^n$ is a polytope with vertex set
$V(P)$ and $B\in\R^n$ is a point.

\begin{theo}[Thales covering theorem]\label{thm:cover}
Let $S\se\R^n$ be a finite nonempty set, $P=\conv S$, and let $B\in\R^n$.
For $V\in S$ put
\[
D_V=\{\,Z\in\R^n:\ \langle Z-V,\,Z-B\rangle\le0\,\},
\]
the closed ball whose diameter is the segment $[V,B]$ (\S\ref{sec:notation}; a
single point if $V=B$). Then
\begin{enumerate}[label=(\alph*)]
\item $P\se\bigcup_{V\in S}D_V$. Equivalently, for every $Z\in P$ there is
      $V\in S$ with $\angle(V,Z,B)\ge\frac{\pi}{2}$, and equivalently with
      \[
      \delta(V,B)^2\ \ge\ \delta(V,Z)^2+\delta(Z,B)^2 .
      \]
\item $P\not\se\bigcup_{V\in S}\inte D_V$; that is, the closed balls
      cannot be replaced by the open ones.
\end{enumerate}
In particular, if $P$ is an $n$-polytope and $S=V(P)$ its vertex set, then
$P$ is covered by the balls whose diameters join a vertex of $P$ to $B$,
and the antipode of $B$ in each ball $D_V$ is the vertex $V$.
\end{theo}
\begin{proo}
(a) Let $Z\in P$ and write $Z=\sum_{V\in S}\lambda_VV$ with $\lambda_V\ge0$ and
$\sum_V\lambda_V=1$. Put $\varphi_V=\langle Z-V,Z-B\rangle$. Since
$V\mapsto\varphi_V$ is an affine function on $\R^n$,
\[
\sum_{V\in S}\lambda_V\,\varphi_V
=\Bigl\langle Z-\sum_{V\in S}\lambda_VV,\ Z-B\Bigr\rangle
=\langle Z-Z,Z-B\rangle=0 .
\]
A weighted average of reals equal to $0$ has a term $\le0$; so $\varphi_V\le0$, i.e.\
$Z\in D_V$, for some $V\in S$.

For the equivalent forms: $\varphi_V\le0$ says that the vectors $V-Z$ and $B-Z$
include an obtuse or right angle, i.e.\ $\angle(V,Z,B)\ge\frac\pi2$ (for
$Z\ne V,B$; for $Z\in\{V,B\}$ we have $\varphi_V=0$).
Moreover, the cosine law for the triangle $V,Z,B$,
\[
\delta(V,B)^2=\delta(V,Z)^2+\delta(Z,B)^2-2\varphi_V ,
\]
turns $\varphi_V\le0$ into $\delta(V,B)^2\ge\delta(V,Z)^2+\delta(Z,B)^2$.

(b) Let $Z\in P$ be a point farthest from $B$. The function
$f=\|\cdot-B\|^2$ is convex with gradient $2(Z-B)$ at $Z$, so for every
$Y\in P$
\[
f(Z)\ \ge\ f(Y)\ \ge\ f(Z)+2\langle Z-B,\,Y-Z\rangle ,
\]
i.e.\ $\langle Z-B,\,Y-Z\rangle\le0$. Taking $Y=V\in S$ gives
$\varphi_V=\langle Z-V,\,Z-B\rangle\ge0$ for all $V\in S$, so $Z$ lies in
none of the open balls $\inte D_V=\{\varphi_V<0\}$.
\end{proo}

\begin{figure}[htbp]
\centering
\begin{tikzpicture}[scale=0.92,line cap=round,line join=round,
                    every node/.style={font=\footnotesize,fill=white,
                                       inner sep=1pt}]
  \coordinate (V1) at (0,0);
  \coordinate (V2) at (4,0);
  \coordinate (V3) at (0,3);
  \coordinate (B)  at (1,1);
  \fill[gray!12] (V1)--(V2)--(V3)--cycle;
  \draw[gray!80,densely dashed,thin] (2.5,0.5) circle (1.58114);
  \draw[gray!80,densely dotted,thin] (0.5,2.0) circle (1.11803);
  \draw[gray!80,solid,thin]          (0.5,0.5) circle (0.70711);
  \draw[gray!45,thin] (V1)--(B) (V2)--(B) (V3)--(B);
  \draw[thick] (V1)--(V2)--(V3)--cycle;
  \gen{V1} \gen{V2} \gen{V3}
  \fill[white] (B) circle (3.6pt); \fill[blue!70!black] (B) circle (2.6pt);
  \node[below left=2pt]  at (V1) {$V_1$};
  \node[below right=2pt] at (V2) {$V_2$};
  \node[left=3pt]        at (V3) {$V_3$};
  \node[right=3pt]       at (B)  {$B$};
  \node[left]  at (-0.25,0.50) {$D_{V_1}$};
  \node[right] at ( 4.12,0.50) {$D_{V_2}$};
  \node[above] at ( 0.50,3.18) {$D_{V_3}$};
  \node[fill=none] at (2,-1.55) {\small (a)};
\end{tikzpicture}\hspace{1.0cm}
\begin{tikzpicture}[scale=0.92,line cap=round,line join=round,
                    every node/.style={font=\footnotesize,fill=white,
                                       inner sep=1pt}]
  \coordinate (V1) at (0,0);
  \coordinate (V2) at (4,0);
  \coordinate (V3) at (0,3);
  \coordinate (B)  at (3,2.5);
  \fill[gray!12] (V1)--(V2)--(V3)--cycle;
  \draw[gray!80,solid,thin]          (1.5,1.25) circle (1.95256);
  \draw[gray!80,densely dashed,thin] (3.5,1.25) circle (1.34629);
  \draw[gray!80,densely dotted,thin] (1.5,2.75) circle (1.52069);
  \draw[gray!45,thin] (V1)--(B) (V2)--(B) (V3)--(B);
  \draw[thick] (V1)--(V2)--(V3)--cycle;
  \gen{V1} \gen{V2} \gen{V3}
  \fill[white] (B) circle (3.6pt); \fill[blue!70!black] (B) circle (2.6pt);
  \node[below left=2pt]  at (V1) {$V_1$};
  \node[below right=2pt] at (V2) {$V_2$};
  \node[left=2pt]        at (V3) {$V_3$};
  \node[above right=3pt] at (B)  {$B$};
  \node[left]  at (-0.47,1.25) {$D_{V_1}$};
  \node[right] at ( 4.87,1.25) {$D_{V_2}$};
  \node[above] at ( 1.50,4.29) {$D_{V_3}$};
  \node[fill=none] at (2,-1.55) {\small (b)};
\end{tikzpicture}
\caption{Theorem~\ref{thm:cover} in the plane, for the triangle
$V_1V_2V_3$ with $V_1=(0,0)$, $V_2=(4,0)$, $V_3=(0,3)$. The closed balls with
diameters $[V_1,B]$, $[V_2,B]$, $[V_3,B]$ (solid, dashed, dotted) all pass through
$B$, and their union contains the whole triangle. (a) $B=(1,1)$ inside the triangle;
the radii are $\tfrac12\sqrt2$, $\tfrac12\sqrt{10}$ and $\tfrac12\sqrt5$. (b)
$B=(3,\tfrac52)$ outside the triangle; the radii are $\tfrac14\sqrt{61}$,
$\tfrac14\sqrt{29}$ and $\tfrac14\sqrt{37}$. In both cases the covering is sharp:
each vertex lies outside the two balls belonging to the other vertices, so no single
ball contains $P$. The figures are exact.}
\label{fig:cover-thales}
\end{figure}
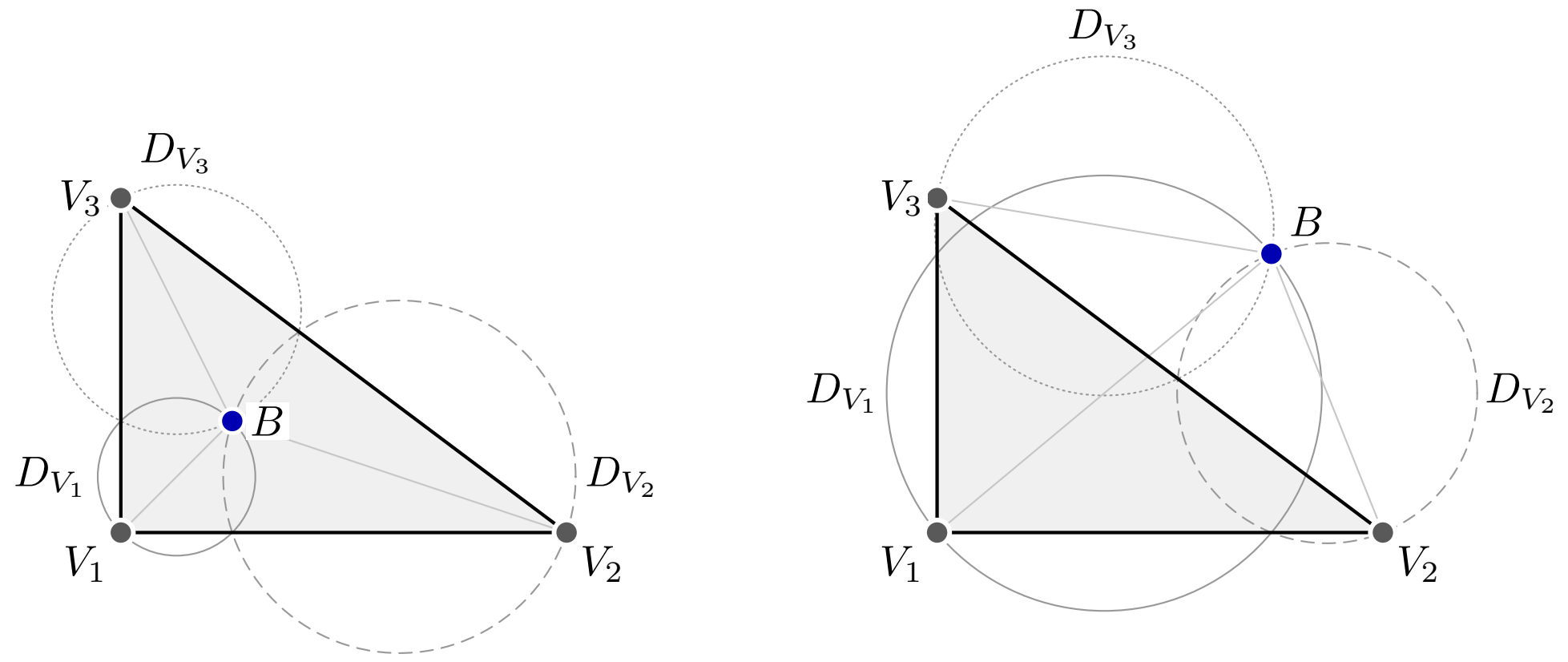

\begin{rem}\label{rem:cover-sharp}
The proof is affine and uses no convexity beyond $P=\conv S$; consequently
$B$ may be any point of $\R^n$, also outside $P$ (Figure~\ref{fig:cover-thales}(b)). The statement has a transparent geometric form:
\[
Z\notin\bigcup_{V\in S}D_V
\quad\Longleftrightarrow\quad
\langle V-Z,\,Z-B\rangle<0\ \text{ for all }V\in S ,
\]
i.e.\ all points of $S$ lie in the open half-space bounded by the
hyperplane through $Z$ perpendicular to $ZB$ and containing $B$. But
$Z\in\conv S$ cannot be strictly separated from $S$ by a hyperplane through
$Z$; that is the whole content of the theorem. Dropping the term
$\delta(V,Z)^2$ in (a) gives $\delta(V,B)\ge\delta(Z,B)$: every point of $P$ has a
vertex at least as far from $B$, whence
$\max_{V\in S}\delta(V,B)=\max_{Z\in P}\delta(Z,B)$. This weaker, distance-only
statement already follows from the convexity of $Z\mapsto\delta(Z,B)$; the theorem
strengthens it to the Pythagorean inequality
$\delta(V,B)^2\ge\delta(V,Z)^2+\delta(Z,B)^2$ for the triangle $V,Z,B$ with a right or
obtuse angle at $Z$.
\end{rem}

The covering of Theorem~\ref{thm:cover} contains a canonical tiling subordinate to it.
Write $\varphi_V(Z)=\langle Z-V,Z-B\rangle$ as in the proof; by \S\ref{sec:notation},
$\varphi_V(Z)=\|Z-M_V\|^2-\rho_V^2$ with $M_V=\frac12(V+B)$ and
$\rho_V=\frac12\delta(V,B)$, i.e.\ $\varphi_V(Z)$ is the \emph{power} of $Z$ with respect
to the Thales sphere $\partial D_V$.

\begin{coro}[power diagram of the Thales spheres]\label{cor:power}
In the situation of Theorem~\ref{thm:cover}, for $V\in S$ let
\[
P_V=\{Z\in P:\ \varphi_V(Z)\le\varphi_W(Z)\ \text{for all }W\in S\}
\]
be the cell of $V$ in the power diagram of the Thales spheres, restricted to $P$. Then
$P_V=P\cap(B+N(V))$, where $N(V)=\{c:\ \langle c,V\rangle\ge\langle c,W\rangle\ \text{for all
}W\in S\}$ is the normal cone of $P$ at $V$; the sets $P_V$ cover $P$ and have pairwise
disjoint interiors; and $P_V\subseteq D_V$ for every $V\in S$.
\end{coro}
\begin{proo}
For $V,W\in S$ we have $\varphi_V(Z)-\varphi_W(Z)=\langle Z-B,\,W-V\rangle$. So
$\varphi_V(Z)\le\varphi_W(Z)$ for all $W$ iff $Z-B\in N(V)$, which gives
$P_V=P\cap(B+N(V))$. The cones $N(V)$ cover $\R^n$ and have pairwise disjoint interiors
(they form the normal fan of $P$; for $V\in S$ not a vertex of $P$ the cone $N(V)$ has
empty interior), which gives the second statement. Finally, for $Z\in P_V$ the proof
of Theorem~\ref{thm:cover}(a) gives
$\varphi_V(Z)=\min_W\varphi_W(Z)\le\sum_W\lambda_W\varphi_W(Z)=0$, i.e.\ $Z\in D_V$.
\end{proo}

\begin{rem}\label{rem:voronoi}
We have $\varphi_V(Z)=\|Z-B\|^2-\langle Z-B,\,V-B\rangle$ and
$\|Z-V\|^2-\|V-B\|^2=\|Z-B\|^2-2\langle Z-B,\,V-B\rangle$. Up to terms independent of
$V$, both decrease with $\langle Z-B,V\rangle$, so the cells $P_V$ are also the cells of
the power diagram of the vertices with the weights $\|V-B\|^2$, i.e.\ of the spheres
centred at the vertices and passing through $B$; all radical hyperplanes of these
spheres pass through their common point $B$. If $P$ is inscribed in a sphere centred at
$B$, the weights are equal and the cells $P_V$ form the \emph{Voronoi tessellation} of $P$
by its vertices.
\end{rem}

For a vertex $V$ of $P$, the orthogonal projections $x_G$ of $B$ onto the faces $G\ni V$ lie
on the sphere bounding $D_V$ (Lemma~\ref{lm:orth-thales}); so the orthogonal cutout at $V$,
when it exists, lies in $D_V$: it is inscribed in the Thales ball of $V$.

\begin{rem}[an interpretation in spatial voting]\label{rem:opinion}
In the spatial model of voting, voters and parties are points of an opinion space
$\R^n$, and Euclidean distance measures the dissimilarity of opinions. Let $W$ be a
voter and $P_1,\dots,P_k$ the parties, $W\neq P_i$. A position $z$ tolerates the party
$P_i$ at least as much as $W$ does if $\delta(z,P_i)\le\delta(W,P_i)$, i.e.\ if $z$ lies in
the ball $B(P_i,\delta(P_i,W))$, whose diameter is $[W,A_i]$ with $A_i=2P_i-W$, the mirror
image of $W$ in $P_i$. Call $z$ a \emph{potential co-voter} of $W$ if this holds for some
party, and $A_i$ the \emph{extreme co-voter} with respect to $P_i$. The set of potential
co-voters is star-shaped with respect to $W$ but in general not convex: for $W=(0,0)$ and
$P_{1,2}=(\mp1,0)$ the two unit discs meet only in $W$, and the midpoint $(0,1)$ of the
boundary points $(\mp1,1)$ lies in neither. Nevertheless, Theorem~\ref{thm:cover} applied to
$S=\{A_1,\dots,A_k\}$ and $B=W$ shows that every convex combination of extreme co-voters,
i.e.\ every compromise among them, is a potential co-voter; by
Theorem~\ref{thm:cover}(b), the closed balls cannot be replaced by open ones. The statement
holds for every subfamily of the parties, for instance for the natural neighbours of $W$,
the parties joined to $W$ by an edge of the Delaunay triangulation of
$\{P_1,\dots,P_k,W\}$. Whether it is of interest in the theory of voting depends on a
comparison with that literature, which we do not attempt here.
\end{rem}

\section{Orthogonal tilings}\label{sec:pedal}

The canonical configuration of \cite[Corollary~6.5]{Bipolar} is orthogonal: its points are
the feet of the perpendiculars from $O$. In this section we consider orthogonal
configurations of an arbitrary interior point, in every dimension. By
Section~\ref{sec:cutouts}, their cutouts are the natural candidates for orthogonally
bipolar cells; we determine when they are such cells and tile the polytope.

In this section we use the orthogonal kernel $\ker_\perp(P)$ of Definition~\ref{def:pedal},
the set of points $B\in\inte P$ whose orthogonal configuration exists
(\cite[Definition~5.7]{Bipolar}). By Proposition~\ref{rem:centres} it is an open convex
polyhedral set, possibly empty, consisting of the points about which $P$ is perfectly
centered. For $B\in\ker_\perp(P)$, Broadie \cite[proof of Lemma~4.6]{Broadie} used the
orthogonal configuration to triangulate the boundary of $P$, in analogy with the
barycentric subdivision; here its points serve as the vertices of the cutouts.

For a vertex $v$ of $P$ and a face $G\ni v$ let $T^v_G=\cone\{w-v:\ w\in G\}$, so that
$v+T^v_G$ is the tangent cone of $G$ at $v$.

For a vertex $v$ let $K_v$ be the set of points $B\in\inte P$ such that
$x_{B,G}\in v+\relint T^v_G$ for every face $G\ni v$; it depends only on the tangent
cone of $P$ at $v$.

\begin{lm}\label{lem:pedal}
$\ker_\perp(P)=\bigcap_{v\in V(P)}K_v$.
\end{lm}
\begin{proo}
A polytope $G$ is the intersection of its tangent cones at its vertices, and so
$\relint G=\bigcap_{v\in G}\bigl(v+\relint T^v_G\bigr)$: every facet of $G$ contains a
vertex $v$, and the corresponding open halfspace of $\aff G$ contains
$v+\relint T^v_G$.
\end{proo}

The lemma does not reduce the number of conditions: both sides are described by
the same linear inequalities, grouped by faces on the left and by vertices on the
right. Its use is that each $K_v$ is a condition on the cutout at $v$ alone. In
dimension three, $B\in K_v$ means exactly that the orthogonal cutout at $v$ is
admissible (Corollary~\ref{cor:lines}); so the orthogonal kernel is the set of points
at which all orthogonal cutouts are admissible. We first show that orthogonally
bipolar cutouts always tile the polytope, in every dimension. Recall (\S\ref{sec:notation})
that this means two things: every point of $P$ lies in some cutout $C_v$, and no point
lies in the interiors of two different cutouts. Both will follow from an explicit
description of the cutouts, which is the key fact about them: each orthogonally
bipolar cutout $C_v$ is the part of $P$ cut out by the normal cone of $P$ at $v$,
translated to $B$ (Lemma~\ref{lem:cell}). The tiling then turns out to be
face-to-face.

\paragraph{Standing assumptions.} For the rest of this section, $P\subset\R^n$ is an
$n$-polytope and $B\in\inte P$; we write $x_G=x_{B,G}$ for the nonempty faces $G$ of $P$.
For a vertex $v$ of $P$, $C_v=\conv\{x_G:\ v\in G\}$ is the \emph{orthogonal cutout} at
$v$, and
\[
N_P(v)=\{c:\ \langle c,v-w\rangle\ge0\ \text{for all neighbours }w\text{ of }v\}
\]
is the normal cone of $P$ at $v$. We say that the cutout $C_v$ is \emph{orthogonally
bipolar} if it is an orthogonally bipolar polytope with poles $v,B$ whose vertices are
the points $x_G$, $v\in G$.

\begin{lm}\label{lem:cell}
Let $B\in\ker_\perp(P)$, and let $v$ be a vertex whose cutout $C_v$ is orthogonally
bipolar. Then $C_v=P\cap\bigl(B+N_P(v)\bigr)$.
\end{lm}
\begin{proo}
Since $B\in\ker_\perp(P)$, each $x_G$ lies in $\relint G$, so $x_G\in K$ iff
$G\subseteq K$, for faces $G,K$ of $P$.

\emph{Faces of $C_v$ at $v$.} Let $G\ni v$ be a proper face of $P$, and let $H$ be a
supporting hyperplane of $P$ with $P\cap H=G$. Since $C_v\subseteq P$ and $v\in C_v\cap H$,
the hyperplane $H$ also supports $C_v$, and $C_v\cap H=C_v\cap G$; so $C_v\cap G$ is a
face of $C_v$ containing $v$ (for $G=P$ this holds trivially); its vertices are the $x_K$ with $v\in K\subseteq G$. These faces are
distinct for distinct $G$. By \cite[Proposition~3.13]{Bipolar}, $y\mapsto[v,y]$ is a bijection
from the vertices of the bipolar polytope $C_v$ onto its faces containing $v$; so $C_v$
has exactly $|\{G:\ v\in G\}|$ faces containing $v$. Hence these are exactly the sets
$C_v\cap G$; in particular $[v,x_G]=C_v\cap G$, and $x_G\le x_K$ in the vertex lattice
of $C_v$ iff $G\subseteq K$. So the atoms of $C_v$ are the points $x_e$, $e\ni v$ an
edge.

\emph{The facets of $C_v$.} Every facet of $C_v$ contains exactly one pole. A facet
containing $v$ is $C_v\cap\Phi$ for a facet $\Phi\ni v$ of $P$, so it lies in $\Phi$. A
facet containing $B$ is $[x_e,B]$ for an edge $e=[v,w]$; by (O3) its affine hull is the
hyperplane through $x_e$ perpendicular to $e$, and since $x_e-B\perp e$ this is
$H_e=\{p:\ \langle p-B,\,v-w\rangle=0\}$. The cell $C_v$ lies on the side of $H_e$
containing $v$, i.e.\ in $\{p:\ \langle p-B,v-w\rangle\ge0\}$. As a polytope is the
intersection of the closed halfspaces of its facets,
\[
C_v=\bigcap_{\Phi\ni v}\Phi^-\ \cap\ \bigcap_{w}\{p:\ \langle p-B,v-w\rangle\ge0\}
=P\cap\bigl(B+N_P(v)\bigr),
\]
where $\Phi^-$ is the halfspace of $P$ bounded by $\aff\Phi$, $w$ runs over the
neighbours of $v$, and $N_P(v)$ is as in the statement. For the second equality note that $P$ lies in every
$\Phi^-$, and that $C_v\subseteq P$.
\end{proo}

\begin{prop}\label{prop:tiling}
Let $B\in\ker_\perp(P)$, and suppose that every cutout $C_v$ is orthogonally bipolar.
Then the cutouts $C_v$, $v\in V(P)$, form a face-to-face tiling of $P$.
\end{prop}
\begin{proo}
By Lemma~\ref{lem:cell}, $C_v=P\cap(B+N_P(v))$ for every vertex $v$.

\emph{Tiling.} The normal cones $N_P(v)$, $v\in V(P)$, form the normal fan of $P$.
First, they cover $\R^n$: for $p\in P$, the vector $p-B$ lies in $N_P(v)$ for a vertex
$v$ maximizing $\langle p-B,\cdot\rangle$ on $P$, so $p\in C_v$. Second, their
interiors are pairwise disjoint, and the interior of $C_v=P\cap(B+N_P(v))$ lies in
$B+\inte N_P(v)$; so no point lies in the interiors of two cells.

\emph{Face-to-face.} Let $v\ne w$ be vertices, let $G$ be the smallest face of $P$
containing both, and fix a point $g\in\relint G$. A vector $c$ lies in
$N_P(v)\cap N_P(w)$ iff the face of $P$ on which $\langle c,\cdot\rangle$ is maximal
contains $v$ and $w$, i.e.\ contains $G$. For $c\in N_P(v)$ this happens iff
$\langle c,v-g\rangle=0$: the inequality $\langle c,v-g\rangle\ge0$ holds since $v$
maximizes $\langle c,\cdot\rangle$, and equality means that $g$ is a maximizer as well,
so that the maximal face contains $g\in\relint G$ and hence $G$. Therefore
\[
C_v\cap C_w=P\cap\bigl(B+N_P(v)\cap N_P(w)\bigr)
=C_v\cap\{p:\ \langle p-B,\,v-g\rangle=0\}.
\]
The linear function $p\mapsto\langle p-B,v-g\rangle$ is nonnegative on $C_v$, because
$p-B\in N_P(v)$ for $p\in C_v$; so $C_v\cap C_w$ is a face of $C_v$. By symmetry, with
$w$ in place of $v$, it is also a face of $C_w$. Hence the tiling is face-to-face.
\end{proo}

\begin{rem}[power diagrams and Voronoi tessellations]\label{rem:power}
By Lemma~\ref{lem:cell}, $C_v=P\cap(B+N_P(v))$, i.e.\ the orthogonally bipolar
cutouts are the cells of the power diagram of the Thales spheres, restricted to $P$
(Corollary~\ref{cor:power}); if $P$ is inscribed in a sphere centred at $B$, they are its
Voronoi cells (Remark~\ref{rem:voronoi}). These cells tile $P$ for every point $B$; the
point of Proposition~\ref{prop:tiling} and Theorem~\ref{thm:pedal} is that for
$B\in\ker_\perp(P)$ they are the orthogonal cutouts. For instance, for a regular polytope and
its centre the orthogonal cutouts are the Voronoi cells of the vertices.

More generally, a tiling of $P$ by convex cells $C_v$, $v\in V(P)$, sharing the vertex $B$
is $P\cap(B+\Sigma)$ for the complete fan $\Sigma$ formed by the tangent cones of the
cells at $B$; it is the restriction of a power diagram with all cells meeting at $B$ iff
$\Sigma$ is the normal fan of a polytope. Indeed, power diagrams
are exactly the regions where one of finitely many affine functions is largest
\cite[Observation~4 and Theorem~4]{Aurenhammer87}, \cite[\S2.2]{Bipolar}; for
cells meeting at $B$ these functions can be taken of the form $\langle p-B,c_v\rangle$,
and the regions where $\langle\cdot,c_v\rangle$ is largest form the normal fan of
$\conv\{c_v\}$. Such tilings are the regular subdivisions in the sense of \cite[\S2.2.3]{DLRS}.
\end{rem}

\begin{theo}[orthogonal tilings]\label{thm:pedal}
The following are equivalent:
\begin{enumerate}[label=(\roman*)]
\item $B\in\ker_\perp(P)$;
\item every cutout $C_v$ is orthogonally bipolar;
\item the cutouts $C_v$ are orthogonally bipolar and form a face-to-face tiling of $P$.
\end{enumerate}
\end{theo}
\begin{proo}
(i)$\Rightarrow$(ii). Let $v$ be a vertex. By (i), $B\in K_v$, i.e.\
$x_{B,G}\in v+\relint T^v_G$ for every face $G\ni v$. These are the hypotheses of
Proposition~\ref{prop:conecut} for the tangent cone of $P$ at $v$, whose faces are the cones
$v+T^v_G$ with $\aff(v+T^v_G)=\aff G$; so the cutout $C_v=\conv\{x_{B,G}:G\ni v\}$ is an
orthogonally bipolar polytope whose vertices are the projections $x_{B,G}$.

(ii)$\Rightarrow$(i). Let $v$ be a vertex and $G\ni v$ a face. The tangent cone of
$C_v$ at $v$ is that of $P$, since the atoms $x_e$ lie on the edges of $P$ at $v$; and
$x_G$ is the vertex of $C_v$ with $\aff[v,x_G]=\aff G$. By
Proposition~\ref{prop:obp-cut}, $x_G$ lies in the relative interior of the face
$v+T^v_G$ of this cone. So $B\in K_v$ for every $v$, and $B\in\ker_\perp(P)$ by
Lemma~\ref{lem:pedal}.

(ii)$\Rightarrow$(iii). By the previous step $B\in\ker_\perp(P)$, so
Proposition~\ref{prop:tiling} applies.

(iii)$\Rightarrow$(ii) is clear.
\end{proo}

The tiling itself is not the point of Theorem~\ref{thm:pedal}: for every $B\in\inte P$
the sets $P\cap(B+N_P(v))$ tile $P$ face-to-face, being the translated normal fan
(Remark~\ref{rem:power}). The content of the theorem is that for $B\in\ker_\perp(P)$ these
cells are exactly the orthogonal cutouts, i.e.\ that the projections $x_G$, $v\in G$,
form the whole vertex set of $P\cap(B+N_P(v))$. The kernel condition only says that
each $x_G$ lies in its own face; a face $G$ could still meet a cone $B+N_P(H)$ with
$H\not\subseteq G$, creating further vertices, and excluding this is the interval
condition (I) for the section of the vertex cone (Theorem~\ref{prop:obpreal}). In
dimension three this follows from \cite[Theorem~5.4(a)]{Bipolar} for polygons; in every
dimension it follows from Broadie's theorem (Remark~\ref{rem:pedal}).

For example, a midscribed $3$-polytope, in particular a canonical one, is perfectly centered
about the centre $O$ of its midsphere \cite[Theorem~5.4(c)]{Bipolar}, so
$O\in\ker_\perp(P)$. By Theorem~\ref{thm:pedal}, the orthogonal cutouts at $O$ then tile $P$
by orthogonally bipolar polytopes $C_v$, each inscribed in its Thales ball $D_v$
(Section~\ref{sec:app}); thus the covering of Theorem~\ref{thm:cover} with $S=V(P)$ and
$B=O$ refines to such a tiling. These cells are those of the canonical configuration of
\cite[Corollary~6.5]{Bipolar}.

\begin{rem}\label{rem:pedal-regular}
For $n=3$ and $B\in\ker_\perp(P)$, the orthogonal configuration is also regular in
the sense of \cite[Definition~6.1]{Bipolar}, so the tiling of Theorem~\ref{thm:pedal} is an
instance of \cite[Theorem~6.3]{Bipolar}. Indeed, (C1) holds by assumption. (C2): for an edge
$e=F\cap F'$ the vectors $x_F-B$ and $x_{F'}-B$ are nonzero multiples of the normals of
$F$ and $F'$, hence independent and orthogonal to $e$, as is $x_e-B$; so
$\Pi_e=\aff\{B,x_F,x_{F'}\}$ is the plane through $B$ perpendicular to $e$, and it
contains $x_e$ but not $e$. (C3): for a vertex $v$ of a facet $F$ with edges $e,e'$ at
$v$, the points $x_e,x_{e'}$ are the feet of the perpendiculars from $x_F$ to the lines
of $e,e'$, so the quadrilateral $v\,x_e\,x_F\,x_{e'}$ has right angles at $x_e,x_{e'}$
and is convex; this is (C3) by \cite[Lemma~6.2]{Bipolar}. (C4): recall that
$W_v=\{x_G:\ v\in G\ne P\}$; for $v\in G\ne P$ we have
$\langle x_G-B,v-B\rangle=\langle x_G-B,v-x_G\rangle+\|x_G-B\|^2=\|x_G-B\|^2>0$, so the
linear function $p\mapsto\langle p-B,v-B\rangle$ is positive on $W_v$ and vanishes at
$B$, whence $B\notin\conv W_v$.
\end{rem}

In particular \cite[Corollary~6.5]{Bipolar} is the case where $P$ is canonical and $B=O$:
the points of tangency and the incentres of the facets lie in the relative interiors.
The orthogonal kernel may be empty, for instance if some facet is a polygon whose
strips perpendicular to the edges have no common point, as in
Figure~\ref{fig:centres}(b); for such polytopes only non-orthogonal configurations
can be regular, as in \cite[Conjecture~6.8]{Bipolar}.

\begin{rem}\label{rem:pedal}
Since $B\in K_v$ depends only on the tangent cone at $v$, Theorem~\ref{thm:pedal} is a
statement about single cones: for the section $Q'$ of the tangent cone by the hyperplane
through $B$ perpendicular to $vB$, the condition $B\in K_v$ says that ($\star1$) holds for
$Q'$ about $B$, and the conclusion is the interval condition for $Q'$ about $B$
(Theorem~\ref{prop:obpreal}). So the implication (i)$\Rightarrow$(ii) for
$n$-polytopes is Broadie's theorem for the $(n-1)$-dimensional sections. Conversely,
Theorem~\ref{thm:pedal} gives Broadie's theorem back: let $Q$ satisfy ($\star1$) about $O$,
and let $P=T\cap\{x_0\le h\}$ be the pyramid of the notation before
Lemma~\ref{lem:tall}, with $h>\beta$ and $B=(\beta,O_d)$. The projections of $B$ onto
the faces of $P$ containing the apex are the points $x_F$, which lie in the relative
interiors of these faces by ($\star1$) and Lemma~\ref{lem:tall}(b), and below the
height $h$; the projections onto the faces of the base $\{h\}\times h\varepsilon Q$ are the
points $(h,h\varepsilon q_G)$, which lie in their relative interiors by ($\star1$). So
$B\in\ker_\perp(P)$, the apex cutout is orthogonally bipolar, and
Theorem~\ref{prop:obpreal} gives the interval condition for $Q$. Thus, in each
dimension, Theorem~\ref{thm:pedal} is equivalent to Broadie's theorem one dimension
lower.
\end{rem}

\section{Open problems}\label{sec:open}

\begin{enumerate}[label=\arabic*.,leftmargin=*]
\item \emph{Orthogonally bipolar realizations.} The face lattice $\F(Q)$ of a polytope
      $Q$ is the vertex lattice of an orthogonally bipolar polytope if and only if the
      combinatorial type of $Q$ has a perfectly centered realization
      Theorem~\ref{prop:obpreal} and Theorem~\ref{conj:I}; this holds in dimension at most
      three, but not for every $4$-polytope \cite[Theorem~3.9]{Dobbins}. Which
      combinatorial types of dimension at least four have a perfectly centered
      realization?
\item \emph{Inscribed versus orthogonally bipolar.} By Theorem~\ref{theo:obpthales}
      every orthogonally bipolar polytope is inscribed in the sphere with diameter
      $\sigma\nu$, and by Remark~\ref{rem:converse} the converse fails in dimension
      three. Describe, for a given combinatorial type, the orthogonally bipolar
      realizations within the space of all geometric bipolar realizations inscribed in
      that sphere; for instance, of which dimension is this subset?
\item \emph{Infinite systems.} Proposition~\ref{prop:bases} applies to every finite
      subfamily of a system of vectors $(b_i)$ in a Hilbert space, and under the sign
      condition each finite subfamily gives a Hjelmslev cuboid whose faces correspond to
      the intervals of a finite Boolean lattice, as for the cube \cite{MetropolisRota}. For
      which infinite systems, for instance Riesz bases, is the closed convex hull
      of the projections of $\nu$ onto the closed spans of all subfamilies an infinite
      Hjelmslev cuboid, i.e.\ a compact convex set whose faces correspond to the intervals
      of the power set of the index set, with all extreme points on the sphere with diameter
      $[0,\nu]$? For orthonormal systems this is the Hilbert cube of
      Remark~\ref{rem:bases}; in general, compactness and the sign condition for infinite
      subfamilies are the issues. In the same spirit: can the iterated interval polytopes of
      Corollary~\ref{cor:iterate} be embedded consistently into $\ell^2$, so that their
      union realizes Lawrence's lattice $\mathcal I_\infty$, which is isomorphic to its
      own lattice of intervals \cite{Lawrence}?
\end{enumerate}

\section*{Acknowledgements}
The author thanks Martin Winter for valuable advice, in particular for pointing out the
connection with perfectly centered polytopes and antiprisms and the work of Broadie,
Fukuda and Weibel, and Dobbins.

\section*{Statements and Declarations}

\noindent\textbf{Funding.} The author acknowledges the support by the Specific Research
Project no.\ FSI-S-26-9581 of Brno University of Technology.

\noindent\textbf{Competing interests.} The author has no competing interests to declare
that are relevant to the content of this article.

\noindent\textbf{Data availability.} No datasets were generated or analysed in this study.
The numerical computations mentioned in the paper only illustrate the results; the
scripts are available from the author on request.

\noindent\textbf{Use of AI tools.} During the preparation of this paper the author used Claude (Anthropic) as an
assistant: for checking and discussing proofs, searching the literature, suggesting
reformulations and restructuring, carrying out numerical experiments, and preparing the
TikZ code of the figures from exact geometric data. All results, proofs, references and
figures were verified by the author, who takes full responsibility for the content.


\begin{thebibliography}{99}

\bibitem{Bipolar} J.~Pavl\'ik, \emph{Face lattices of polytopes and
bipolar polytopes}, preprint, arXiv:2611.13896, 2026.

\bibitem{Carr2005} R.~D.~Carr, G.~Konjevod, Polyhedral combinatorics, in:
H.~J.~Greenberg (ed.), \emph{Tutorials on Emerging Methodologies and Applications in
Operations Research}, Int.\ Ser.\ Oper.\ Res.\ Manag.\ Sci.~76, Springer, New York,
2005, pp.~2-1--2-46.

\bibitem{Schrijver95} A.~Schrijver, Polyhedral combinatorics, in:
\emph{Handbook of Combinatorics}, Vol.~2 (R.~L.~Graham, M.~Gr\"otschel,
L.~Lov\'asz, eds.), Elsevier, Amsterdam, 1995, pp.~1649--1704.

\bibitem{MetropolisRota} N.~Metropolis, G.-C.~Rota, Combinatorial structure
of the faces of the $n$-cube, \emph{SIAM J.\ Appl.\ Math.}\ 35 (1978),
689--694.

\bibitem{Grunbaum} B.~Gr\"unbaum, \emph{Convex Polytopes}, 2nd ed., Springer
GTM 221, 2003.

\bibitem{Ziegler} G.~M.~Ziegler, \emph{Lectures on Polytopes}, Springer GTM
152, 1995.

\bibitem{McMullenSchulte} P.~McMullen, E.~Schulte, \emph{Abstract Regular
Polytopes}, Cambridge Univ.\ Press, 2002.

\bibitem{Aurenhammer87} F.~Aurenhammer, Power diagrams: properties, algorithms and
  applications, \emph{SIAM J.~Comput.} 16 (1987), 78--96.
\bibitem{Aurenhammer91} F.~Aurenhammer, Voronoi diagrams---a survey of a fundamental
  geometric data structure, \emph{ACM Comput.~Surv.} 23 (1991), 345--405.
\bibitem{DLRS} J.~A.~De~Loera, J.~Rambau, F.~Santos, \emph{Triangulations:
  Structures for Algorithms and Applications}, Algorithms and Computation in
  Mathematics~25, Springer, Berlin, 2010.
\bibitem{Broadie} M.~Broadie, A theorem about antiprisms, \emph{Linear Algebra Appl.}~66
  (1985), 99--111.
\bibitem{FukudaWeibel} K.~Fukuda, C.~Weibel, $f$-vectors of Minkowski additions of
  convex polytopes, \emph{Discrete Comput.\ Geom.}~37 (2007), 503--516.
\bibitem{Lawrence} J.~Lawrence, Dual-antiprisms and partitions of powers of~$2$ into
  powers of~$2$, \emph{Discrete Comput.\ Geom.}~61 (2019), 465--478.
\bibitem{Dobbins} M.~G.~Dobbins, Antiprismlessness, or: reducing combinatorial
  equivalence to projective equivalence in realizability problems for polytopes,
  \emph{Discrete Comput.\ Geom.}~57 (2017), 966--984.
\bibitem{GrunbaumAnti} B.~Gr\"unbaum, Are prisms and antiprisms really boring?
  (Part~3), \emph{Geombinatorics}~15 (2005), no.~2, 69--78.
\bibitem{GrunbaumAnti2} B.~Gr\"unbaum, Are prisms and antiprisms really boring?
  (Part~2), \emph{Geombinatorics}~14 (2004), 10--17.
\bibitem{BermanPlemmons} A.~Berman, R.~J.~Plemmons, \emph{Nonnegative Matrices in the
  Mathematical Sciences}, Classics in Applied Mathematics~9, SIAM, Philadelphia, 1994.
\end{thebibliography}
\end{document}